\documentclass[12pt]{article}
\usepackage{amsmath,amssymb,amsthm,latexsym,amscd}

\title{Algebras of distributions
of binary isolating formulas of a complete theory\footnote{\indent \ \ \ The work is supported by RFBR (grant
12-01-00460-a).}}
\author{Ilya V. Shulepov\footnote{salvodore@mail.ru} \ and Sergey V.
Sudoplatov\footnote{sudoplat@math.nsc.ru}}
\date{May 15, 2012}
\begin{document}
\maketitle

\begin{abstract}
We define a class of algebras describing links of binary isolating
formulas on a set of realizations for a family of $1$-types of a
complete theory. We prove that a set of labels for binary
isolating formulas on a set of realizations for a $1$-type $p$
forms a groupoid of a special form if there is an atomic model
over a realization of $p$. We describe the class of these
groupoids and consider features of these groupoids in a general
case and for special theories. A description of the class of
partial groupoids relative to families of $1$-types is given.

\end{abstract}

\bigskip
In \cite{SuLP} (see also \cite{Su041}--\cite{Su08}), a series of
constructions is introduced admitting to realize key properties of
countable theories and to obtain a classification of countable
models of small (in particular, of Ehrenfeucht) theories with
respect to two basic characteristics: Rudin--Keisler preorders and
distribution functions for numbers of limit models. The
construction of these theories is essentially based on the
definition of special directed graphs with colored vertices and
arcs as well as on the definition of $(n+1)$-ary predicates that
turn prime models over realizations of $n$-types to prime models
over realizations of $1$-types and reducing links between prime
models over finite sets to links between prime models over
elements such that these links are defined by principal arcs and
edges.

In the paper, we consider a general approach to the description of
binary links between realizations of $1$-types in terms of labels
of pairwise non-equivalent isolating formulas, being represented
implicitly or for some special cases in \cite{SuLP}--\cite{Shul3}.
This approach is naturally interpretable in the class of relation
partial algebras \cite{HiHo, Mad}.

In Section 1, we define a class of algebras distributing binary
isolating formulas and introduce preliminary definitions,
notations, and properties of algebras connected with relations of
isolation and semi-isolation. In Sections 2, we describe some
basic examples for these algebras and for types basing these
algebras. In Section 3, we define a groupoid
$\mathfrak{P}_{\nu(p)}$ of principal formulas on a set of
realizations of $1$-type $p$ (assuming that there is an atomic
model over a realization of $p$) with respect to a {\em regular}
labelling function $\nu(p)$ for pairwise non-equivalent principal
formulas $\varphi(a,y)$ for which $\varphi(a,x)\vdash p(x)$ holds,
$\models p(a)$. In Section 4, we collect the basic properties of
groupoids $\mathfrak{P}_{\nu(p)}$ and the significant subgroupoids
of $\mathfrak{P}_{\nu(p)}$. In Section 5, using the
successively-annihilating sums we construct two kinds of monoids
$\mathfrak{P}_{\nu(p)}$ containing an arbitrary group. In Section
6, we produce a list of properties characterizing the class of
groupoids $\mathfrak{P}_{\nu(p)}$. Features of these groupoids for
the class of special theories are exposed in Section 7. In Section
8, we define the notion of join of groupoids and show the
mechanism of extension of basic properties of
$\mathfrak{P}_{\nu(p)}$ to the class of partial groupoids being
joins of groupoids $\mathfrak{P}_{\nu(p)}$. In final Section 8, we
produce a list of properties characterizing the class of partial
groupoids correspondent to algebras of distributions for binary
isolating formulas on a family of types.

We use the standard relation algebraic, model-theoretical,
semigroup, and graph-theoretic terminology \cite{HiHo}--\cite{SO2}
as well as some notions, notations, and constructions in
\cite{SuLP}.

\medskip
\section{Preliminary notions, notations \\ and properties}

\medskip
{\bf Definition} \cite{SuLP, Su08, BSV}. Let $T$ be a complete
theory, $\mathcal{M}\models T$. Consider types $p(x),q(y)\in
S(\varnothing)$, realized in $\mathcal{M}$, and all {\em
$(p,q)$-preserving}\index{Formula!$(p,q)$-preserving} formulas
$\varphi(x,y)$ of $T$, i.~e., formulas for which there is $a\in M$
such that $\models p(a)$ and $\varphi(a,y)\vdash q(y)$. Now, for
each such a formula $\varphi(x,y)$, we define a binary relation
$R_{p,\varphi,q}\rightleftharpoons\{(a,b)\mid\mathcal{M}\models
p(a)\wedge\varphi(a,b)\}.$ If $(a,b)\in R_{p,\varphi,q}$, $(a,b)$
is called a {\em $(p,\varphi,q)$-arc}\index{$(p,\varphi,q)$-arc}.
If $\varphi(a,y)$ is principal (over $a$), the $(p,\varphi,q)$-arc
$(a,b)$ is also {\em
principal}\index{$(p,\varphi,q)$-arc!principal}\index{Arc!principal}.
If,~in~addition, $\varphi(x,b)$ is principal (over $b$), the set
$[a,b]\rightleftharpoons\{(a,b),(b,a)\}$ is said to be a {\em
principal $(p,\varphi,q)$-edge}\index{$(p,\varphi,q)$-edge
principal}\index{Edge!principal}. $(p,\varphi,q)$-arcs and
$(p,\varphi,q)$-edges are called {\em arcs}\index{Arc} and {\em
edges}\index{Edge} respectively if we say about fixed or some
formula $\varphi(x,y)$. If~$(a,b)$ is a principal arc and $(b,a)$
is not a principal arc (on any formula) then $(a,b)$ is called
{\em irreversible}\index{Arc!principal!irreversible}.

\medskip
For types $p(x),q(y)\in S(\varnothing)$, we denote by ${\rm
PF}(p,q)$\index{${\rm PF}(p,q)$} the set
$$\{\varphi(x,y)\mid\varphi(a,y)\mbox{
is a principal formula, } \varphi(a,y)\vdash q(y),\mbox{ where
}\models p(a)\}.$$ Let ${\rm PE}(p,q)$\index{${\rm PE}(p,q)$} be
the set of pairs of formulas $(\varphi(x,y),\psi(x,y))\in{\rm
PF}(p,q)$ such that for any (some) realization $a$ of $p$ the sets
of solutions for $\varphi(a,y)$ and $\psi(a,y)$ coincide. Clearly,
${\rm PE}(p,q)$ is an equivalence relation on the set ${\rm
PF}(p,q)$. Notice that each ${\rm PE}(p,q)$-class $E$ corresponds
to either a principal edge or to an irreversible principal arc
connecting realizations of $p$ and $q$ by some (any) formula in
$E$. Thus the quotient ${\rm PF}(p,q)/{\rm PE}(p,q)$ is
represented as a disjoint union of sets ${\rm
PFS}(p,q)$\index{${\rm PFS}(p,q)$} and ${\rm
PFN}(p,q)$\index{${\rm PFN}(p,q)$}, where ${\rm PFS}(p,q)$
consists of ${\rm PE}(p,q)$-classes correspondent to principal
edges and ${\rm PFN}(p,q)$ consists of ${\rm PE}(p,q)$-classes
correspondent to irreversible principal arcs.

The sets ${\rm PF}(p,p)$, ${\rm PE}(p,p)$, ${\rm PFS}(p,p)$ , and
${\rm PFN}(p,p)$ are denoted by ${\rm PF}(p)$,\index{${\rm
PF}(p)$} ${\rm PE}(p)$,\index{${\rm PE}(p)$} ${\rm
PFS}(p)$,\index{${\rm PFS}(p)$} and ${\rm PFN}(p)$\index{${\rm
PFN}(p)$} respectively.

Let $T$ be a complete theory without finite models,
$U=U^-\,\dot{\cup}\,\{0\}\,\dot{\cup}\,U^+$ be an alphabet of
cardinality $\geq|S(T)|$ and consisting of {\em negative
elements}\index{Element!negative} $u^-\in U^-$,\index{$U^-$} {\em
positive elements}\index{Element!positive} $u^+\in
U^+$\index{$U^+$}, and zero $0$. As usual, we write $u<0$ for any
$u\in U^-$ and $u>0$ for any $u\in U^+$.\footnote{If $U$ is at
most countable, we assume that $U$ is a subset of the set $\mathbb
Z$ of integers.} The set $U^-\cup\{0\}$ is denoted by $U^{\leq
0}$\index{$U^{\leq 0}$} and $U^+\cup\{0\}$ is denoted by $U^{\geq
0}$\index{$U^{\geq 0}$}. Elements of $U$ are called {\em
labels}.\index{Label}

Let \ $\nu(p,q)\mbox{\rm : }{\rm PF}(p,q)/{\rm PE}(p,q)\to U$ \ be
\ injective \ {\em labelling \
functions},\index{Function!labelling} $p(x),q(y)\in
S(\varnothing)$, for which negative elements correspond to the
classes in ${\rm PFN}(p,q)/{\rm PE}(p,q)$ and non-negative
elements correspond to the classes in ${\rm PFS}(p,q)/{\rm
PE}(p,q)$ such that $0$ is defined only for $p=q$ and is
represented by the formula $(x\approx y)$,
$\nu(p)\rightleftharpoons\nu(p,p)$\index{$\nu(p)$}. We
additionally assume that $\rho_{\nu(p)}\cap\rho_{\nu(q)}=\{0\}$
for $p\ne q$ (where, as usual, we denote by $\rho_f$ the image of
the function $f$) and
$\rho_{\nu(p,q)}\cap\rho_{\nu(p',q')}=\varnothing$ if $p\ne q$ and
$(p,q)\ne(p',q')$. Labelling functions with the properties above
as well families of these functions are said to be {\em
regular}\index{Function!labelling!regular}\index{Family!regular}.
Further we shall consider only regular labelling functions and
their regular families.

We denote by $\theta_{p,u,q}(x,y)$\index{$\theta_{p,u,q}(x,y)$}
formulas in ${\rm PF}(p,q)$ with a label $u\in\rho_{\nu(p,q)}$. If
the type $p$ is fixed and $p=q$ then the formula
$\theta_{p,u,q}(x,y)$ is denoted by
$\theta_u(x,y)$.\index{$\theta_u(x,y)$}

Note that if $\theta_{p,u,q}(x,y)$ and $\theta_{q,v,p}(x,y)$ are
formulas witnessing that for realizations $a$ and $b$ of $p$ and
$q$ respectively the pairs $(a,b)$ and $(b,a)$ are principal arcs
then the formula $\theta_{p,u,q}(x,y)\wedge\theta_{q,v,p}(y,x)$
witnesses that $[a,b]$ is a principal edge. Moreover the
(non-negative) label $v$ corresponds uniquely to the {\em
invertible}\index{Label!invertible} label $u$ and vice versa. The
labels $u$ and $v$ are {\em reciprocally
inverse}\index{Labels!reciprocally inverse} and are denoted by
$v^{-1}$\index{$v^{-1}$} and $u^{-1}$\index{$u^{-1}$}
respectively.

For types $p_1,p_2,\ldots,p_{k+1}\in S^1(\varnothing)$ and sets
$X_1,X_2,\ldots,X_k\subseteq U$ of labels we denote
by\index{$P(p_1,X_1,p_2,X_2,\ldots,p_k,X_k,p_{k+1})$}
$$P(p_1,X_1,p_2,X_2, \ldots, p_k,X_k,p_{k+1})$$ the set
of all labels $u\in U$ correspondent to formulas
$\theta_{p_1,u,p_{k+1}}(x,y)$ satisfying, for realizations $a$ of
$p_1$ and some $u_1\in X_1,\ldots,u_k\in X_k$, the following
condition:
$$
\theta_{p_1,u,p_{k+1}}(a,y)\vdash\theta_{p_1,u_1,p_2,u_2,\ldots,p_k,u_k,p_{k+1}}(a,y),$$
where\index{$\theta_{p_1,u_1,p_2,u_2,\ldots,p_k,u_k,p_{k+1}}(x,y)$}
$$\theta_{p_1,u_1,p_2,u_2,\ldots,p_k,u_k,p_{k+1}}(x,y)\rightleftharpoons$$
$$\rightleftharpoons\exists x_2,x_3,\ldots
x_{k-1}(\theta_{p_1,u_1,p_2}(x,x_2)\wedge\theta_{p_2,u_2,p_3}(x_2,x_3)\wedge\ldots$$
$$\ldots\wedge\theta_{p_{k-1},u_{k-1},p_k}(x_{k-1},x_k)\wedge\theta_{p_k,u_k,p_{k+1}}(x_k,y)).
$$

Thus the Boolean $\mathcal{P}(U)$ of $U$ is the universe of an
{\em algebra of distributions of binary isolating
formulas}\index{Algebra!of distributions of binary isolating
formulas} with $k$-ary operations
$$P(p_1,\cdot,p_2,\cdot, \ldots, p_k,\cdot,p_{k+1}),$$
where $p_1,\ldots,p_{k+1}\in S^1(\varnothing)$. This algebra has a
natural restriction to any family $R\subseteq S^1(\varnothing)$.

Note that if some set $X_i$ is disjoint with
$\rho_{\nu(p_i,p_{i+1})}$, in particular, if it is empty then
$$P(p_1,X_1,p_2,X_2,\ldots,p_k,X_k,p_{k+1})=\varnothing.$$

Note also that if $X_i\not\subseteq \rho_{\nu(p_i,p_{i+1})}$ for
some $i$ then
$$P(p_1,X_1,p_2,X_2, \ldots, p_k,X_k,p_{k+1})=$$
$$=P(p_1,X_1\cap\rho_{\nu(p_1,p_2)},p_2,X_2\cap\rho_{\nu(p_1,p_2)},
\ldots, p_k,X_k\cap\rho_{\nu(p_k,p_{k+1})},p_{k+1}).$$

In view of the previous equation, it is enough to assume
$X_i\subseteq \rho_{\nu(p_i,p_{i+1})}$, $i=1,\ldots,k$ for the
values $P(p_1,X_1,p_2,X_2,\ldots,p_k,X_k,p_{k+1})$.

If each set $X_i$ is a singleton consisting of an element $u_i$
then we use $u_i$ instead of $X_i$ in
$P(p_1,X_1,p_2,X_2,\ldots,p_k,X_k,p_{k+1})$ and
write\index{$P(p_1,u_1,p_2,u_2,\ldots,p_k,u_k,p_{k+1})$}
$$P(p_1,u_1,p_2,u_2,\ldots,p_k,u_k,p_{k+1}).$$

By the definition the following equality holds:
$$P(p_1,X_1,p_2,X_2, \ldots, p_k,X_k,p_{k+1})=$$
$$=\cup\{P(p_1,u_1,p_2,u_2, \ldots, p_k,u_k,p_{k+1})\mid u_1\in X_1,\ldots, u_k\in
X_k\}.$$ Hence the specification of $P(p_1,X_1,p_2,X_2, \ldots,
p_k,X_k,p_{k+1})$ is reduced to the specifications of
$P(p_1,u_1,p_2,u_2, \ldots, p_k,u_k,p_{k+1})$. Note also that
$P(p,X,q)=X$ for any $X\subseteq\rho_{\nu(p,q)}$.

Clearly, if $u_i=0$ then $p_i=p_{i+1}$ for nonempty sets
$$P(p_1,u_1,p_2,u_2,\ldots,p_i,0,p_{i+1},\ldots,p_k,u_k,p_{k+1})$$
and the following conditions hold:
$$
P(p_1,0,p_1)=\{0\},
$$
$$P(p_1,u_1,p_2,u_2,\ldots,p_i,0,p_{i+1},\ldots, p_k,u_k,p_{k+1})=$$
$$=P(p_1,u_1,p_2,u_2,\ldots,p_i,u_{i+1},p_{i+2},\ldots,
p_k,u_k,p_{k+1}).$$

If all types $p_i$ equal to a type $p$ then we write
$P_p(X_1,X_2,\ldots,X_k)$\index{$P_p(X_1,X_2,\ldots,X_k)$} and
$P_p(u_1,u_2,\ldots,u_k)$\index{$P_p(u_1,u_2,\ldots,u_k)$} as well
as $\lfloor X_1,X_2,\ldots,X_k\rfloor _p$\index{$\lfloor
X_1,X_2,\ldots,X_k\rfloor _p$} and $\lfloor
u_1,u_2,\ldots,u_k\rfloor _p$\index{$\lfloor
u_1,u_2,\ldots,u_k\rfloor _p$} instead of
$$P(p_1,X_1,p_2,X_2,\ldots,p_k,X_k,p_{k+1})$$ and
$$P(p_1,u_1,p_2,u_2,\ldots,p_k,u_k,p_{k+1})$$
respectively.  We omit the index ${\cdot}_p$ if the type $p$ is
fixed. In this case, we write
$\theta_{u_1,u_2,\ldots,u_k}(x,y)$\index{$\theta_{u_1,u_2,\ldots,u_k}(x,y)$}
instead of $\theta_{p,u_1,p,u_2,\ldots,p,u_k,p}(x,y)$.

\medskip
{\bf Definition} (A.~Pillay\index{Pillay A.} \cite{Pi83}). Let
$\mathcal{M}$ be a model of a theory~$T$, $\bar{a}$ and $\bar{b}$
be tuples in $\mathcal{M}$, $A$ be a subset of $M$. The tuple
$\bar{a}$ \emph{semi-isolates}\index{Tuple!semi-isolates a tuple}
the tuple $\bar{b}$ over the set $A$ if there exists a formula
$\varphi(\bar{a},\bar{y})\in{\rm tp}(\bar{b}/A\bar{a})$ for which
$\varphi(\bar{a},\bar{y})\vdash{\rm tp}(\bar{b}/A)$ holds. In this
case we say that the formula $\varphi(\bar{a},\bar{y})$ (with
parameters in $A$) \emph{witnesses that $\bar{b}$ is semi-isolated
over $\bar{a}$ with respect to $A$}.\index{Formula!witnessing!on
semi-isolation}

Similarly, a tuple $\bar{a}$ \emph{isolates}\index{Tuple!isolates
a tuple} a tuple $\bar{b}$ over $A$ if there exists a formula
$\varphi(\bar{a},\bar{y})\in{\rm tp}(\bar{b}/A\bar{a})$ for which
$\varphi(\bar{a},\bar{y})\vdash{\rm tp}(\bar{b}/A)$ and
$\varphi(\bar{a},\bar{y})$ is a principal (i.~e., isolating)
formula. In this case we say that the formula
$\varphi(\bar{a},\bar{y})$ (with parameters in $A$)
\emph{witnesses that $\bar{b}$ is isolated over $\bar{a}$ with
respect to $A$}.\index{Formula!witnessing!on isolation}

If $\bar{a}$ (semi-)isolates $\bar{b}$ over $\varnothing$, we
simply say that $\bar{a}$ {\em {\rm (}semi-{\rm
)}isolates}\index{Tuple!semi-isolates a
tuple}\index{Tuple!isolates a tuple} $\bar{b}$; and if a formula
$\varphi(\bar{a},\bar{y})$ witnesses that $\bar{a}$
(semi-)isolates $\bar{b}$ over $\varnothing$ then we say that
$\varphi(\bar{a},\bar{y})$ \emph{witnesses that $\bar{a}$ {\rm
(}semi-{\rm )}isolates $\bar{b}$}.\index{Formula!witnessing!on
semi-isolation}\index{Formula!witnessing!on isolation}

If $q\in S(T)$ then ${\rm SI}_q$\index{${\rm SI}_q$} (in the model
$\mathcal{M}$) denotes the relation of semi-isolation (over
$\varnothing$) on a set of realizations of $q$:
$$
{\rm
SI}_q\rightleftharpoons\{(\bar{a},\bar{b})\mid\mathcal{M}\models
q(\bar{a})\wedge q(\bar{b}) \mbox{ and }\bar{a}\mbox{
semi-isolates }\bar{b}\}.
$$

Similarly, we denote by $I_q$\index{$I_q$} (in the model
$\mathcal{M}$) the relation of isolation (over $\varnothing$) on a
set of realizations of $q$:
$$
I_p\rightleftharpoons\{(\bar{a},\bar{b})\mid\:\mathcal{M}\models
q(\bar{a})\wedge q(\bar{b}) \mbox{ and }\bar{a}\mbox{ isolates
}\bar{b}\}.
$$

For a family $R\subset S(T)$ of $1$-types we denote by
$I_R$\index{$I_R$} (in the model $\mathcal{M}$) the set
$$\{(a,b)\mid {\rm tp}(a),{\rm tp}(b)\in R\mbox{ and }a\mbox{
isolates }b\}$$ and by ${\rm SI}_R$\index{${\rm SI}_R$} (in
$\mathcal{M}$) the set
$$\{(a,b)\mid {\rm tp}(a),{\rm tp}(b)\in R\mbox{ and }a\mbox{
semi-isolates }b\}.$$

Clearly, $I_R\subseteq{\rm SI}_R$ and, for any set of realizations
of types in $R$, the relations $I_R$ and ${\rm SI}_R$ are
reflexive. As shown in \cite{Pi83}, the relation of semi-isolation
on the set of tuples in an arbitrary model is transitive and, in
particular, any relation ${\rm SI}_R$ is transitive.

\medskip
{\bf Lemma 1.1} \cite{Ki}--\cite{Cas11}. {\em $(1)$ If a tuple
$\bar{a}$ isolates a tuple $\bar{b}$, whereas $\bar{b}$ does not
isolate $\bar{a}$, then $\bar{b}$ does not semi-isolate $\bar{a}$.

$(2)$ If $(a,b)\in I_R$ and $(b,a)\in{\rm SI}_R$ then $(b,a)\in
I_R$.}

\medskip
{\em Proof.} (1) Suppose that $\varphi(\bar{a},\bar{y})$ isolates
${\rm tp}(\bar{b}/\bar{a})$. Assume the contrary (i.~e., $\bar{b}$
semi-isolates $\bar{a}$) and take a formula
$\psi(\bar{x},\bar{b})$ witnessing that $\bar{b}$ semi-isolates
$\bar{a}$. Now as ${\rm tp}(\bar{a}/\bar{b})$ is nonisolated,
there exists a formula $\chi(\bar{x},\bar{y})$ such that
$\varphi(\bar{x},\bar{b})\wedge\psi(\bar{x},\bar{b})\wedge\chi(\bar{x},\bar{b})$
and
$\varphi(\bar{x},\bar{b})\wedge\psi(\bar{x},\bar{b})\wedge\neg\chi(\bar{x},\bar{b})$
are both consistent. Moreover both formulas imply ${\rm
tp}(\bar{a})$. Hence
$\varphi(\bar{a},\bar{y})\wedge\chi(\bar{a},\bar{y})$ and
$\varphi(\bar{a},\bar{y})\wedge\neg\chi(\bar{a},\bar{y})$ are both
consistent. This contradicts the fact that
$\varphi(\bar{a},\bar{y})$ is a principal formula.

(2) follows immediately from (1).~$\Box$

\medskip
{\bf Proposition 1.2.} {\em $(1)$ If $p,q\in R$ are principal
types then $\rho_{\nu(p,q)}\cup\rho_{\nu(q,p)}\subseteq U^{\geq
0}$.

$(2)$ If $p,q\in R$, $p$ is a principal type and $q$ is a
non-principal type then $\rho_{\nu(p,q)}=\varnothing$ and
$\rho_{\nu(q,p)}\subseteq U^-$.}

\medskip
{\em Proof.} (1) If $\rho_{\nu(p,q)}$ contains a label $u<0$ then
there are realizations $a$ and $b$ of $p$ and $q$ respectively
such that $(a,b)\in I_R$ and $(b,a)\notin I_R$. So by Lemma 1.1,
$(b,a)\notin{\rm SI}_R$. But since $p(x)$ contains a principal
formula $\varphi(x)$, this formula witnesses that $(b,a)\in{\rm
SI}_R$. The contradiction implies that $\rho_{\nu(p,q)}\subseteq
U^{\geq 0}$. Similarly we obtain $\rho_{\nu(q,p)}\subseteq U^{\geq
0}$.

(2) Let $\varphi(x)$ be a principal formula of $p(x)$. If $\models
p(a)$, $\models q(b)$, and $(a,b)\in I_R$ that witnessed by a
formula $\theta_u(x,y)$, the formula $\exists
x(\varphi(x)\wedge\theta_u(x,y))$ isolates $q(y)$. Since $q$ is
not isolated we obtain $\rho_{\nu(p,q)}=\varnothing$. By the same
reason, $\rho_{\nu(q,p)}\subseteq U^-$.~$\Box$

\medskip
{\bf Corollary 1.3.} {\em If $p(x)$ is a principal type then
$\rho_{\nu(p)}\subseteq U^{\geq 0}$.}

\medskip
{\bf Proposition 1.4.} {\em Let $p_1,p_2,\ldots,p_{k+1}$ be types
in $S^1(\varnothing)$. The following assertions hold.

$(1)$ If $u_i\in\rho_{\nu(p_i,p_{i+1})}$, $i=1,\ldots,k$, and some
$u_i$ is negative then
$$P(p_1,u_1,p_2,u_2, \ldots, p_k,u_k,p_{k+1})\subseteq U^-.$$

$(2)$ If $u_i\in\rho_{\nu(p_i,p_{i+1})}$, $i=1,\ldots,k$, and all
elements $u_i$ are not negative then
$$P(p_1,u_1,p_2,u_2, \ldots,
p_k,u_k,p_{k+1})\subseteq U^{\geq 0}.$$

$(3)$ If $u_i\in\rho_{\nu(p_i,p_{i+1})}$, $i=1,\ldots,k$, and all
elements $u_i$ are non-negative, then all elements of the set
$$X\rightleftharpoons P(p_1,u_1,p_2,u_2,\ldots,p_k,u_k,p_{k+1})$$
are invertible and the set $X^{-1}\rightleftharpoons\{v^{-1}\mid
v\in X\}$\index{$X^{-1}$} coincides with the set
$P(p_{k+1},u^{-1}_k,p_k,u^{-1}_{k-1},\ldots,p_2,u^{-1}_1,p_{1})$.}

{\em Proof.} (1) Let $v$ be a label in $P(p_1,u_1,p_2,u_2, \ldots,
p_k,u_k,p_{k+1}).$ Consider realizations $a_i$ of $p_i$ such that
$$\models\theta_{p_i,u_i,p_{i+1}}(a_i,a_{i+1}),\,\,\,i=1,\ldots,
k,\,\,\,\models\theta_{p_1,v,p_{k+1}}(a_1,a_{k+1}).$$ For the
family $R=\{p_1,p_2,\ldots,p_{k+1}\}$ we have $(a_1,a_{k+1})\in
I_R$, $(a_i,a_{i+1})\in I_R$, $i=1,\ldots,k$, and so
$(a_i,a_j)\in{\rm SI}_R$ for $i\leq j$. If $u_i<0$ then
$(a_{i+1},a_i)\notin I_R$ and then, by Lemma 1.1,
$(a_{i+1},a_i)\notin{\rm SI}_R$. If $v\geq 0$ then
$(a_{k+1},a_1)\in I_R$ and, by transitivity of ${\rm SI}_R$ and
$(a_{i+1},a_{k+1}),(a_{k+1},a_1),(a_1,a_i)\in{\rm SI}_R$ we get
$(a_{i+1},a_i)\in{\rm SI}_R$ that is impossible. Since the element
$v\in P(p_1,u_1,p_2,u_2, \ldots, p_k,u_k,p_{k+1})$ is taken
arbitrarily the set $P(p_1,u_1,p_2,u_2, \ldots, p_k,u_k,p_{k+1})$
consists of negative elements.

(2) Take again elements $v,a_1,a_2,\ldots,a_{k+1}$ as for (1). If
$u_i\geq 0$ then $(a_{i+1},a_i)\in I_R$, $i=1,\ldots,k$. By
transitivity of the relation ${\rm SI}_R$, the element $a_{k+1}$
semi-isolates the element $a_1$. In view of $(a_1,a_{k+1})\in
I_R$, by Lemma 1.1, we have $(a_{k+1},a_1)\in I_R$ and so $v\geq
0$. Since the element $$v\in P(p_1,u_1,p_2,u_2, \ldots,
p_k,u_k,p_{k+1})$$ is taken arbitrarily the set
$P(p_1,u_1,p_2,u_2, \ldots, p_k,u_k,p_{k+1})$ consists of
non-negative elements.

(3) follows immediately from (2).~$\Box$

\medskip
{\bf Corollary 1.5.} {\em Restrictions of $U$ to the sets $U^{\leq
0}$ and $U^{\geq 0}$ form subalgebras of the algebra of
distributions of binary isolating formulas. Each element of the
restriction to $U^{\geq 0}$ has a unique inverse element. The
operation of inversion is coordinated with the operations of the
algebra.}

\medskip
\section{Examples}

\medskip
Consider some examples for distributions of labels of binary
isolating formulas on sets of realizations of types $p(x)\in
S(\varnothing)$ for countable theories $T$.

I. If $|\rho_{\nu(p)}|=1$ then $(x\approx y)$ is the unique
principal formula up to equivalence. It is possible only in the
following cases:

(1) $T$ is small (i.~e., with countable $S(\varnothing)$) and
satisfies some of the following condition:

(a)  $p(x)$ is a principal type with the only realization;

(b) $p(x)$ is a non-principal type such that if a set
$\{\varphi(a,y)\wedge\neg(a\approx y)\}\cup p(y)$ is consistent,
where $\varphi(x,y)$ is a formula of $T$, $\models p(a)$, then
$\varphi(a,y)\not\vdash p(y)$;

(2) $T$ is a theory with continuum many types and for any formula
$\varphi(x,y)$ of $T$ and for a realization $a$ of $p(x)$ if the
set $\{\varphi(a,y)\wedge\neg(a\approx y)\}\cup p(y)$ is
consistent and $\varphi(a,y)\vdash p(y)$ then there are no
isolating formulas $\psi(a,y)$ such that
$\psi(a,y)\vdash\varphi(a,y)\wedge\neg(a\approx y)$.

The case 1,a is represented by a type being realized by a
constant; the cases 1,b and 2 are represented by theories of unary
predicates with non-principal types $p(x)$ and having countably
many and continuum many types respectively.

\medskip
II. Let $\rho_{\nu(p)}=\{0,1\}$. Then $1^{-1}=1$ and any
realization $a$ of $p$ is linked with the only realization $b$ of
$p$ for which $\models\theta_1(a,b)$ and, moreover,
$\models\theta_1(b,a)$. Then the set of realizations of $p$ splits
on two-element equivalence classes consisting of $\theta_1$-edges.
If $p$ is a principal type of a small theory then a
$\theta_1$-edge is unique, and if $p$ is non-principal the number
of this edges can vary from $1$ to the infinity depending on a
model of a theory.

\medskip
III. Let $\rho_{\nu(p)}=\{-1,0\}$ be a set for a small theory $T$.
By Corollary 1.3 the type $p(x)$ is non-principal and the formula
$\theta_{-1}(x,y)$ witnesses that ${\rm SI}_p$ is non-symmetric.
The formula $\theta_{-1,-1}(x,y)\rightleftharpoons\exists
z(\theta_{-1}(x,z)\wedge\theta_{-1}(z,y))$ is also witnessing that
${\rm SI}_p$ is non-symmetric. By assumption the formula
$\theta_{-1,-1}(a,y)$ is equivalent to the formula
$\theta_{-1}(a,y)$. It means that, on a set of realizations of
$p$, the relation described by the formula
$\theta_{-1}(x,y)\vee(x\approx y)$ is an infinite partial order.
This partial order is dense since if the element $a$ has a
covering element then the formula $\theta_{-1}(a,y)$ is equivalent
to the disjunction of consistent formulas
$\theta_{-1}(a,y)\wedge\theta_{-1,-1}(a,y)$ and
$\theta_{-1}(a,y)\wedge\neg\theta_{-1,-1}(a,y)$, but it is
impossible for the principal formula $\theta_{-1}(a,y)$.

We consider, as a theory with $\rho_{\nu(p)}=\{-1,0\}$, the
Ehrenfeucht's theory $T$, i.~e. the theory of a structure
$\mathcal{M}$, formed from the structure $\langle\mathbb
Q;<\rangle$ by adding constants $c_k$, $c_k<c_{k+1}$,
$k\in\omega$, such that $\lim\limits_{k\to\infty}c_k=\infty$. The
type $p(x)$, isolated by the set of formulas $c_k<x$,
$k\in\omega$, has exactly two non-equivalent isolating formulas:
$\theta_{-1}(a,y)=(a<y)$ and $\theta_0(a,y)=(a\approx y)$, where
$\models p(a)$.

\medskip
IV. Let $\rho_{\nu(p)}=\{-1,0,1\}$. Realizing this equation, we
consider the Ehrenfeucht's example, where each element $a$ is
replaced by an $<$-antichain consisting of two elements $a'$ and
$a''$ such that $\models\theta_1(a',a'')\wedge\theta_1(a'',a')$.
Then we have the following equations for the type $p(x)$ isolated
by the set of formulas $c'_k<x$, $k\in\omega$:
$P_p(-1,-1)=P_p(-1,1)=P_p(1,-1)=\{-1\}$, $P_p(1,1)=\{0\}$.

\medskip
V. The equation $\rho_{\nu(p)}=\{-2,-1,0\}$ with
$P_p(-2,-2)=\{-2\}$ and
$$P_p(-2,-1)=P_p(-1,-2)=P_p(-1,-1)=\{-1\}$$ can be fulfilled by
two dense strict orders $<_1$ and $<_2$ on a set of realizations
of a non-principal type such that $<_1$ immerses $<_2$: $<_1\circ
<_2\,\,=\,\,<_2\circ <_1\,\,=\,\, <_1$.

\medskip
VI. Consider a dense linearly ordered set
$\mathcal{M}=\langle\emph{Q},<\rangle$, $T={\rm Th}(\mathcal{M})$,
and  the unique $1$-type $p$ of $T$. Define a labelling function
$\nu(p)$, for which $0$ corresponds to the formula $(x\approx y)$,
$1$ to $(x < y)$, and $2$ to $(y < x)$. We have $\rho_{\nu(p)}=
\{0,1,2\}$, $P_p(1,2)=P_p(2,1)=\rho_{\nu(p)}$, $P_p(1,1)=\{1\}$,
$P_p(2,2)=\{2\}$.

\medskip
VII. Take a group $\langle G;\,\ast\rangle$ and define, on the set
$G$ binary predicates $Q_g$, $g\in G$, by the following rule:
$$
Q_g=\{(a,b)\in G^2\mid a\ast g=b\}.
$$
If $p(x)$ is a type (of a theory $T$) realized in any model
$\mathcal{M}\models T$ containing $G$ exactly by elements in $G$
connected by definable relations $Q_g$, then the type $p$ is
isolated, the set $G$ is finite, and $\rho_{\nu(p)}$ consists of
non-negative elements bijective with elements in $G$. If
$\rho_{\nu(p)}$ consists of non-negative elements, is bijective
with $G$, and the set of realizations of a principal type $p$ is
not fixed, then, assuming the smallness of the theory, the set $G$
is infinite and the number of connected components with respect to
the relation $Q\rightleftharpoons\bigcup\limits_{g\in G}Q_g$ is
not bounded. At last if the type $p$ is not isolated then the
number of $Q$-components on sets of realizations of $p$ is also
unbounded although the set $G$ can be finite.

The Cayley table of the group $\langle G;\,\ast\rangle$ defines
operations $P_p(\cdot,\ldots,\cdot)$ on the set $\rho_{\nu(p)}$ in
accordance with links between the relations $Q_g$.

\medskip
VIII. Applying to a concrete group we consider the structure
$\mathcal{M}\rightleftharpoons\langle\mathbb Z;s^{(1)}\rangle$
with the unary \emph{successor function}\index{Function!successor}
$s\mbox{\rm : }\mathbb Z\leftrightarrow\mathbb Z$, where
$s(n)=n+1$ for each $n\in\mathbb Z$. For the unique $1$-type $p$
of the theory ${\rm Th}(\mathcal{M})$ the set of pairwise
non-equivalent formulas $\theta_u(x,y)$ is exhausted by the list:
$y\approx\underbrace{s\ldots s}_{n\mbox{\scriptsize \ times}}(x)$
and $x\approx\underbrace{s\ldots s}_{n\mbox{\scriptsize \
times}}(y)$, $n\in\omega$. The set $\rho_{\nu(p)}$ consists of
non-negative elements linked by additive group of integers.

\medskip
\section{Algebra of distributions of binary isolating
formulas on a set of realizations of a type}
\medskip

We consider a complete theory $T$, a type $p(x)\in S(T)$, a
regular labelling function \ $\nu(p)\mbox{\rm : }{\rm PF}(p)/{\rm
PE}(p)\to U$, \ and a family of sets $P_p(u_1,\ldots,u_k)$,
$u_1,\ldots,u_k\in\rho_{\nu(p)}$, $k\in\omega$, of labels for
binary isolating formulas.

Further we denote by $\mathcal{M}_p$ and by $\mathcal{M}(a)$ an
atomic model over a realization $a$ of $p$.

Below we prove some basic properties for sets
$$\lfloor u_1,\ldots,u_k\rfloor \rightleftharpoons P_p(u_1,\ldots,u_k).$$

\medskip
{\bf Proposition 3.1.} {\em {\rm 1.} A set $\lfloor
u_1,u_2\rfloor$ is nonempty if and only if for a realization $a$
of $p$ and for some formula $\theta_v(x,y)$,
$\theta_v(a,y)\vdash\theta_{u_1,u_2}(a,y)$ holds.

{\rm 2.} If a model $\mathcal{M}_p$ exists then the set $\lfloor
u_1,u_2\rfloor$ is nonempty for any $u_1,u_2\in\rho_{\nu(p)}$.

{\rm 3.} The set $\lfloor u_1,u_2,u_3\rfloor$  is nonempty if and
only if for a realization $a$ of $p$ and for some formula
$\theta_v(x,y)$, $\theta_v(a,y)\vdash\theta_{u_1,u_2,u_3}(a,y)$
holds.

{\rm 4.} For any $u_1,u_2,u_3\in\rho_{\nu(p)}$ the following
inclusions are satisfied: $$\lfloor\lfloor
u_1,u_2\rfloor,u_3\rfloor\subseteq\lfloor u_1,u_2,u_3\rfloor,$$
$$\lfloor
u_1,\lfloor u_2,u_3\rfloor\rfloor\subseteq\lfloor
u_1,u_2,u_3\rfloor.$$

{\rm 5.} For any $u_1,u_2,u_3\in\rho_{\nu(p)}$ the inclusion
$$\lfloor
u_1,u_2,u_3\rfloor\subseteq\lfloor\lfloor
u_1,u_2\rfloor,u_3\rfloor$$ holds if and only if for any $v\in
\lfloor u_1,u_2,u_3\rfloor$ there is $v'\in\lfloor u_1,u_2\rfloor$
such that $v\in \lfloor v',u_3\rfloor$.

{\rm 6.} {\rm (Left
semi-associativity)}\index{Semi-associativity!left} If a model
$\mathcal{M}_p$ exists then, for any
$u_1,u_2,u_3\in\rho_{\nu(p)}$,
$$\lfloor\lfloor u_1,u_2\rfloor,u_3\rfloor=\lfloor
u_1,u_2,u_3\rfloor.$$

{\rm 7.} For any $u_1,u_2,u_3\in\rho_{\nu(p)}$ the inclusion
$$\lfloor
u_1,u_2,u_3\rfloor\subseteq\lfloor u_1,\lfloor
u_2,u_3\rfloor\rfloor$$ is true if and only if for any $v\in
\lfloor u_1,u_2,u_3\rfloor$ there is $v'\in\lfloor u_2,u_3\rfloor$
such that $v\in \lfloor u_1,v'\rfloor$.

{\rm 8.} {\rm (Criterion for right
semi-associativity)}\index{Semi-associativity!right} If the model
$\mathcal{M}(a)$ exists, where $\models p(a)$, then for any
$u_1,u_2,u_3\in\rho_{\nu(p)}$ the equality
$$\lfloor u_1,\lfloor
u_2,u_3\rfloor\rfloor=\lfloor u_1,u_2,u_3\rfloor$$ holds if and
only if for any $v\in\lfloor u_1,u_2,u_3\rfloor$ the formula
$\theta_{u_1}(a,y_1)\wedge\theta_{u_2,u_3}(y_1,y_2)\wedge\theta_v(a,y_2)$
is realized in $\mathcal{M}(a)$ by a principal arc $(b_1,b_2)$.

{\rm 9.} {\rm ($(\geq 0)$-associativity)}\index{$(\geq
0)$-associativity} If the model $\mathcal{M}(a)$ exists, where
$\models p(a)$, then for any $u_1,u_2,u_3\in\rho_{\nu(p)}$, where
$u_1\geq 0$,
$$\lfloor \lfloor u_1,u_2\rfloor ,u_3\rfloor =\lfloor u_1,u_2,u_3\rfloor =\lfloor u_1,\lfloor u_2,u_3\rfloor \rfloor.$$}

\medskip
{\em Proof.} 1, 2, 3, 5, 7 follow immediately by the definition.
In view of 4, 8 is an easy reformulation of 7.

4. For the proof of $\lfloor\lfloor
u_1,u_2\rfloor,u_3\rfloor\subseteq\lfloor u_1,u_2,u_3\rfloor$, we
take an arbitrary element $v\in \lfloor \lfloor u_1,u_2\rfloor
,u_3\rfloor $. Then $v\in \lfloor v',u_3\rfloor$ for some $v'\in
\lfloor u_1,u_2\rfloor $, and for any realization $a$ of $p$ we
have
\begin{equation} \theta_{v'}(a,x_2)\vdash \theta_{u_1,u_2}(a,x_2),
\end{equation}
\begin{equation} \theta_{v}(a,y)\vdash\theta_{v',u_3}(a,y).
\end{equation}
By (1), we obtain
\begin{equation} \theta_{v',u_3}(a,y)\vdash
\theta_{u_1,u_2,u_3}(a,y).
\end{equation}
Thus, (2) and (3) imply
$$\theta_{v}(a,y)\vdash\theta_{u_1,u_2,u_3}(a,y),$$
and, consequently, $v\in\lfloor u_1,u_2,u_3\rfloor$.

Now we prove the inclusion $\lfloor u_1,\lfloor
u_2,u_3\rfloor\rfloor\subseteq\lfloor u_1,u_2,u_3\rfloor$. Take an
arbitrary element $v\in\lfloor u_1,\lfloor u_2,u_3\rfloor\rfloor$.
Then $v\in \lfloor u_1,v'\rfloor$ for some $v'\in \lfloor
u_2,u_3\rfloor $, and for any realization $a$ of $p$ we have
\begin{equation} \theta_{v'}(a,y)\vdash \theta_{u_2,u_3}(a,y),
\end{equation}
\begin{equation} \theta_{v}(a,y)\vdash\theta_{u_1,v'}(a,y).
\end{equation}
By (4), we obtain
\begin{equation}
\theta_{u_1,v'}(a,y)\vdash \theta_{u_1,u_2,u_3}(a,y).
\end{equation}
Thus, (5) and (6) imply
$$\theta_{v}(a,y)\vdash\theta_{u_1,u_2,u_3}(a,y),$$
and, consequently, $v\in\lfloor u_1,u_2,u_3\rfloor$.

6. Take a realization $a$ of $p$ and an element $v\in \lfloor
u_1,u_2,u_3\rfloor$. Then, for the principal formula
$\theta_v(a,y)$, we have
$\theta_v(a,y)\vdash\theta_{u_1,u_2,u_3}(a,y)$ and so
$$\mathcal{M}(a)\models\theta_{u_1}(a,b_1)\wedge\theta_{u_2}(b_1,b_2)\wedge\theta_{u_3}(b_2,c)\wedge\theta_v(a,c)$$ для
for some realizations $b_1,$ $b_2$, and $c$ of $p$. Since the
model $\mathcal{M}(a)$ is atomic over $a$ we have
$\theta_{v'}(a,x_2)\vdash\theta_{u_1,u_2}(a,x_2)$ and
$\mathcal{M}(a)\models\theta_{v'}(a,b_2)$ for some $v'\in\lfloor
u_1,u_2\rfloor$. Then $\theta_v(a,y)\vdash\theta_{v',u_3}(a,y)$
and hence $v\in\lfloor v',u_3\rfloor$. Since the element
$v\in\lfloor u_1,u_2,u_3\rfloor $ is chosen arbitrarily, we
obtain, by 5, $\lfloor u_1,u_2,u_3\rfloor \subseteq\lfloor \lfloor
u_1,u_2\rfloor ,u_3\rfloor$ that implies, by 4, the equality
$\lfloor \lfloor u_1,u_2\rfloor ,u_3\rfloor=\lfloor
u_1,u_2,u_3\rfloor$.

9. By 4 and 6, it suffices to prove $\lfloor
u_1,u_2,u_3\rfloor\subseteq\lfloor u_1,\lfloor u_2,u_3\rfloor
\rfloor$ for any $u_1,u_2,u_3\in\rho_{\nu(p)}$, where $u_1\geq 0$.
Let $v$ be an arbitrary element in $\lfloor u_1,u_2,u_3\rfloor$.
Since $u_1\geq 0$ there is the label $u_1^{-1}$ and, in
$\mathcal{M}(a)$, there are realizations $b,c,d$ of $p$ such that
$$
\mathcal{M}(a)\models\theta_{u_1^{-1}}(a,b)\wedge\theta_{u_2}(a,c)\wedge\theta_{u_2,u_3}(a,d)\wedge\theta_v(b,d).
$$
Since the type ${\rm tp}(d/a)$ is principal, we have
$\mathcal{M}(a)\models\theta_{v'}(a,d)$ for some label $v'$. As
$v'\in\lfloor u_2,u_3\rfloor$ and $v\in\lfloor u_1,v'\rfloor$ we
obtain $v\in\lfloor u_1,\lfloor u_2,u_3\rfloor \rfloor$.

If $u_2\geq 0$ and $u_3\geq 0$ we also have the required inclusion
by the following arguments. Since, by Proposition 1.4 (2), $v\geq
0$, and there is a non-negative element $v^{-1}\in\lfloor
u^{-1}_3,u^{-1}_2,u^{-1}_1\rfloor$, then, by 6, we have
$v^{-1}\in\lfloor\lfloor
u^{-1}_3,u^{-1}_2\rfloor,u^{-1}_1\rfloor$. Applying Proposition
1.4 (3), we obtain $v\in\lfloor u_1,\lfloor u_2,u_3\rfloor
\rfloor$.~$\Box$

\medskip
Proposition 3.1 implies

\medskip
{\bf Corollary 3.2.} {\em If there is a model $\mathcal{M}(a)$,
where $\models p(a)$, then the following conditions hold:

{\rm 1.} For any $u_1,u_2,u_3\in\rho_{\nu(p)}$, the equalities
$$\lfloor \lfloor u_1,u_2\rfloor ,u_3\rfloor =\lfloor u_1,u_2,u_3\rfloor\supseteq\lfloor u_1,\lfloor u_2,u_3\rfloor \rfloor$$
are satisfied.

{\rm 2. (Criterion of associativity)} For any
$u_1,u_2,u_3\in\rho_{\nu(p)}$, the equality
$$\lfloor \lfloor u_1,u_2\rfloor ,u_3\rfloor=\lfloor u_1,\lfloor u_2,u_3\rfloor \rfloor$$
hold if and only if $u_1\geq 0$ or, for any $v\in\lfloor
u_1,u_2,u_3\rfloor$, the formula
$\theta_{u_1}(a,y_1)\wedge\theta_{u_2,u_3}(y_1,y_2)\wedge\theta_v(a,y_2)$
is realized in $\mathcal{M}(a)$ by a principal arc $(b_1,b_2)$.}

\medskip
Note that if $\mathcal{M}_p$ does not exist the associativity (as
well as semi-associativ\-ities) can be failed. For instance, if
$\lfloor u_1,u_2\rfloor=\varnothing$ then $\lfloor \lfloor
u_1,u_2\rfloor ,u_3\rfloor$ is also empty although $\lfloor
u_1,\lfloor u_2,u_3\rfloor \rfloor\ne\varnothing$ is admissible.

By Proposition 3.1, having $\mathcal{M}_p$ the associativity can
be failed only by some labels $u_1,u_2,u_3$ with $u_1<0$. By
Proposition 1.4 (1), in this case any label $v\in \lfloor
u_1,u_2,u_3\rfloor$ is also negative. The mechanism presented in
the following example shows that the fault of right
semi-associativity is admitted for any distribution of signs for
nonzero labels $u_2,u_3$: there are small theories with
\begin{equation}
\lfloor u_1,u_2,u_3\rfloor\ne\lfloor u_1,\lfloor
u_2,u_3\rfloor\rfloor.
\end{equation}

\medskip
{\bf Example 3.1.} Obtaining (7) with $u_1<0$, $u_2,u_3\ne 0$, and
a label $v\in\lfloor u_1,u_2,u_3\rfloor\setminus \lfloor
u_1,\lfloor u_2,u_3\rfloor\rfloor $ (i.~e., by Proposition 3.1, 8,
for the non-realizability of the formula
$\varphi(a,y_1,y_2)\rightleftharpoons\theta_{u_1}(a,y_1)\wedge\theta_{u_2,u_3}(y_1,y_2)\wedge\theta_v(a,y_2)$
by principal arcs) we consider the schema of the realization of a
non-$p$-principal $(2,p)$-type in a model $\mathcal{M}_p$ of small
theory presented in \cite[Example 1.3.1]{SuLP} (see also
\cite{Su92}). Defining the type $p(x)$ we introduce a $Q_{u_1}$-
and {\em $Q_v$-ordered} (for binary predicates $Q_{u_1}$ and $Q_v$
correspondent to the labels $u_1$ and $v$) coloring ${\rm Col
}\mbox{\rm : }M_0\to\omega\cup\{\infty\}$ of some graph $\Gamma$
producing unary predicates ${\rm Col}_n=\{a\in M_0\mid {\rm
Col}(a)=n\}$, $n\in\omega$, such that:

(a) for any $m\leq n<\omega$ there are elements $a,b\in M_0$ for
which $\models{\rm Col }_m(a)\wedge{\rm Col }_n(b)\wedge Q(a,b)$;

(b) if $m<n<\omega$ then there are no elements $c,d\in M_0$ for
which $\models{\rm Col}_m(c)\wedge{\rm Col }_n(d)\wedge Q(d,c)$.

Moreover, using a generic construction for $\Gamma$ we obtain the
unique non-principal $1$-type $p(x)$ and it is isolated by the set
$\{\neg{\rm Col}_n(x)\mid n\in\omega\}$.

For each label $u_i$, $i\in\{2,3\}$, depending on its label, we
define a binary predicate $Q_{u_i}$ linking only the same elements
in color if $u_i$ is positive, and with the $Q_{u_i}$-ordering of
${\rm Col}$ if $u_i<0$. Now we introduce labels $v'_n$,
$n\in\omega$, being negative if $u_2<0$ or $u_3<0$ and positive
otherwise, such that $\lfloor u_2,u_3\rfloor=\{v'_n\mid
n\in\omega\}$. We define pairwise disjoint predicates $Q_{v'_n}$
linking only the same elements in color if $v'_n>0$, and linking
with the $Q_{v'_n}$-ordering of ${\rm Col}$ if $v'_n<0$. Moreover,
we require the following condition: for any element $a_k$ of color
$k$ the formula $\varphi(a_k,y_1,y_2)$ is realized by principal
$Q_{v'_n}$-arcs exactly with $n\geq k$. It means that, for
$\models p(a)$, the formula $\varphi(a,y_1,y_2)$ is not realized
by principal arcs, since this formula witnesses that the
non-$p$-principal $(2,p)$-type
$$q(y_1,y_2)\rightleftharpoons p(y_1)\cup
p(y_2)\cup\{\theta_{u_2,u_3}(y_1,y_2)\}\cup\{\neg\theta_{v'_n}(y_1,y_2)\mid
n\in\omega\}$$ is realized in $\mathcal{M}_p$.~$\Box$

\medskip
If the model $\mathcal{M}_p$ exists then, using the left
semi-associativity, by induction on the number of brackets one
prove that all operations $\lfloor\cdot,\cdot,\ldots,\cdot\rfloor$
acting on sets in
$\mathcal{P}(\rho_{\nu(p)})\setminus\{\varnothing\}$ are generated
by the binary operation $\lfloor \cdot,\cdot\rfloor $ on the set
$\mathcal{P}(\rho_{\nu(p)})\setminus\{\varnothing\}$ If we have
the right semi-associativity, the values $\lfloor
X_1,X_2,\ldots,X_k\rfloor$,
$X_1,X_2,\ldots,X_k\subseteq\rho_{\nu(p)}$, do not depend on
sequences of placements of brackets for
$$X_{i,i+1,\ldots,i+m+n}\rightleftharpoons\lfloor
X_{i,i+1,\ldots,i+m},X_{i+m+1,i+m+2,\ldots,i+m+n}\rfloor,$$ where
$X_{1,2,\ldots,k}=\lfloor X_1,X_2,\ldots,X_k\rfloor $.

Thus, having $\mathcal{M}_p$, the groupoid $\mathfrak
P_{\nu(p)}\rightleftharpoons\langle\mathcal{P}(\rho_{\nu(p)})\setminus\{\varnothing\};\lfloor
\cdot,\cdot\rfloor \rangle$\index{$\mathfrak P_{\nu(p)}$}, being a
{\rm (}left{\rm )} {\em
semi-associative}\index{Algebra!semi-associative} algebra, admits
to represent all operations $\lfloor
\cdot,\cdot,\ldots,\cdot\rfloor$ by terms of the language $\lfloor
\cdot,\cdot\rfloor$. Further the operation $\lfloor
\cdot,\cdot\rfloor$ will be also denoted by $\cdot$ and we shall
write $uv$ instead of $u\cdot v$. If the right semi-associativity
fails we shall assume, for $u_1u_2\ldots u_k$, the following
distribution of parentheses: $(((u_1\cdot u_2)\cdot\ldots)\cdot
u_k)$.

Since by the choice of the label $0$ for the formula $(x\approx
y)$ the equalities $X\cdot\{0\}=X$ and $\{0\}\cdot X=X$ are true
for any $X\subseteq\rho_{\nu(p)}$, the groupoid $\mathfrak
P_{\nu(p)}$ has the unit $\{0\}$, and it is a monoid if the
algebra is right semi-associative. We have
$$
Y\cdot Z=\bigcup\{yz\mid y\in Y,z\in Z\}
$$
for any sets
$Y,Z\in\mathcal{P}(\rho_{\nu(p)})\setminus\{\varnothing\}$ in this
structure.

Thus the following proposition holds.

\medskip
{\bf Proposition 3.3.} {\em For any complete theory $T$, any type
$p\in S(T)$ having the model $\mathcal{M}_p$, and the regular
labelling function $\nu(p)$, any operation
$P_p(\cdot,\cdot,\ldots,\cdot)$ on the set
$\mathcal{P}(\rho_{\nu(p)})\setminus\{\varnothing\}$ is
interpretable by a term of the groupoid $\mathfrak P_{\nu(p)}$.}

\medskip
The groupoid $\mathfrak P_{\nu(p)}$ is called the {\em groupoid of
binary isolating formulas over the labelling function
$\nu(p)$}\index{Groupoid!of binary isolating formulas} or the {\em
$I_{\nu(p)}$-groupoid}\index{$I_{\nu(p)}$-groupoid}.

\medskip
Propositions 1.4 and 3.1 imply

\medskip
{\bf Proposition 3.4.} {\em  For any complete theory $T$, any type
$p\in S(T)$ having the model $\mathcal{M}_p$, and the regular
labelling function $\nu(p)$, the restriction of the groupoid
$\mathfrak P_{\nu(p)}$ to the set of non-positive {\rm
(}respectively non-negative{\rm )} labels is a semi-associative
subalgebra of $\mathfrak P_{\nu(p)}$ with the unit $\{0\}$ {\rm
(}and, moreover, it is a monoid{\rm )}.}

\medskip
\section{Characterization of transitivity \\ for the
relation $I_p$. Deterministic, \\ almost deterministic
$I_{\nu(p)}$-groupoids \\ and elements}
\medskip

The following assertion gives a characterization of transitivity
of the relation $I_p$. For simplicity we formulate and prove it
for a $1$-type $p$ although the proof implies the validity for any
complete type $r$ of a theory with a model $\mathcal{M}_r$.

\medskip
{\bf Proposition 4.1.} {\em Let $p(x)$ be a complete type of
complete theory $T$ having a model $\mathcal{M}_p$, $\nu(p)$ be a
regular labelling function. The following conditions are
equivalent:

$(1)$ the relation $I_p$ {\rm (}on a set of realizations of $p$ in
a model $\mathcal{M}\models T${\rm )} is transitive;

$(2)$ for any labels $u_1,u_2\in\rho_{\nu(p)}$ the set
$P_p(u_1,u_2)$ is finite.}

\medskip
{\em Proof.} Let $a,b,c$ be realizations of $p$ such that
$(a,b)\in I_p$ and $(b,c)\in I_p$ witnessed by isolating formulas
$\theta_{u_1}(a,y)$ and $\theta_{u_2}(b,y)$. If the set
$P_p(u_1,u_2)$ is finite and consists of labels $v_1,\ldots,v_k$
then, by existence of $\mathcal{M}_p$, the formula
$\theta_{u_1,u_2}(a,y)$ is equivalent to the formula
$\bigvee\limits_{i=1}^k\theta_{v_i}(a,y)$. Since
$\models\theta_{u_1,u_2}(a,c)$ we have
$\models\bigvee\limits_{i=1}^k\theta_{v_i}(a,c)$ and hence
$\models\theta_{v_i}(a,c)$ for some $i$. Thus, $(a,c)\in I_p$ and
it is witnessed by the formula $\theta_{v_i}(x,y)$. In view of
arbitrary choice of elements $a,b,c$ the implication
$(2)\Rightarrow(1)$ is true.

Now, we assume that, for some $u_1,u_2\in\rho_{\nu(p)}$, the set
$P_p(u_1,u_2)$ is infinite. Then by compactness, for a realization
$a$ of $p$, the set
$$q(a,y)\rightleftharpoons\{\theta_{u_1,u_2}(a,y)\}\cup\{\neg\theta_v(a,y)\mid
v\in P_p(u_1,u_2)\}$$ is consistent. Consider realizations $b$ and
$c$ of $p$ such that
$\models\theta_{u_1}(a,b)\wedge\theta_{u_2}(b,c)$ and $\models
q(a,c)$. We have $(a,b)\in I_p$, $(b,c)\in I_p$, and $(a,c)\notin
I_p$ by the construction of $q$. Thus the relation $I_p$ is not
transitive and we obtain $(1)\Rightarrow(2)$.~$\Box$

\medskip
{\bf Definition.} A structure $\mathfrak P_{\nu(p)}$ is called
({\em almost}) {\em
deterministic},\index{Structure!deterministic}\index{Structure!almost
deterministic}\index{Monoid!deterministic}\index{Monoid!almost
deterministic} if the set $\lfloor u_1,u_2\rfloor$ is a singleton
(is nonempty and finite) for any $u_1,u_2\in\rho_{\nu(p)}$.

\medskip
{\bf Proposition 4.2.} {\em If there is a model $\mathcal{M}_p$
and the structure $\mathfrak P_{\nu(p)}$ is almost deterministic
then $\mathfrak P_{\nu(p)}$ is a monoid.}

\medskip
{\em Proof.} As noticed in Proposition 3.1, the unique obstacle,
for $\mathfrak P_{\nu(p)}$ to be a monoid, can be only the
existence of labels $u_1,u_2,u_3,v$, $u_1<0$, $v<0$, for which
$v\in\lfloor u_1,u_2,u_3\rfloor$ and there are no $v'\in\lfloor
u_2,u_3\rfloor$ with $v\in\lfloor u_1,v'\rfloor$. But, by the
hypothesis, the set $\lfloor u_2,u_3\rfloor$ consists of finitely
many labels $v_1,\ldots,v_k$. Now we take in $\mathcal{M}(a)$,
where $\models p(a)$, elements $b,c,d$ such that
$$
\mathcal{M}(a)\models\theta_{u_1}(a,b)\wedge\theta_{u_2}(b,c)\wedge\theta_{u_3}(c,d)\wedge\theta_v(a,d).
$$
Since the formula $\theta_{u_2,u_3}(b,y)$ is equivalent to the
formula $\bigvee\limits_{i=1}^k\theta_{v_i}(b,y)$, there is a
required label $v'=v_i$ such that
$\mathcal{M}(a)\models\theta_{v'}(b,d)$.~$\Box$

\medskip
{\bf Example 4.1.} \ By \ the \ definition \ any \
polygonometrical \ theory \ ${\rm Th}({\rm
pm}(G_1,G_2,\mathcal{P}))$ (see \cite{SuGP}) has a unique $1$-type
$p(x)\in S(\varnothing)$ and, thus, the structure $\mathfrak
P_{\nu(p)}$ is a monoid with non-negative labels. The (almost)
determinacy of $\mathfrak P_{\nu(p)}$ means that the group $G_1$
of sides is unit or the group $G_2$ of angles is unit
(finite).~$\Box$

\medskip
Any deterministic structure $\mathfrak P_{\nu(p)}$ is a monoid
(being almost deterministic). It is generated by the monoid
$\mathfrak P'_{\nu(p)}=\langle \rho_{\nu(p)};\,
\odot\rangle$,\index{$\mathfrak P'_{\nu(p)}$} where $\lfloor
u,v\rfloor=\{u\odot v\}$ for $u,v\in \rho_{\nu(p)}$.

Thus, the deterministic monoids can be defined by usual Cayley
tables for monoids on a set of labels in $U$ while the almost
deterministic monoids are represented by one-to-finite functions
with two arguments, i.~e., by ternary predicates with finitely
many third coordinates for fixed first and second coordinates.

Considering deterministic structures $\mathfrak P$, being
restrictions of the monoid $\mathfrak P_{\nu(p)}$ to some
subalphabets $U_0$ of the alphabet $U$, we denote by $\mathfrak
P'$ the {\em generating}\index{Monoid!generating} monoid $\langle
U_0;\, \odot\rangle$ such that $\lfloor u,v\rfloor\cap
U_0=\{u\odot v\}$ for $u,v\in U_0$.

The following proposition is a reformulation of Proposition 4.1.

\medskip
{\bf Proposition 4.3.} {\em Let $p(x)$ be a complete type of a
theory $T$ having a model $\mathcal{M}_p$, $\nu(p)$ be a regular
labelling function. The following conditions are equivalent:

$(1)$ the relation $I_p$ {\rm (}on a set of realizations of $p$ in
a model $\mathcal{M}\models T${\rm )} is transitive;

$(2)$ the structure $\mathfrak P_{\nu(p)}$ is an almost
deterministic monoid.}

\medskip
Note that there are no principal edges linking distinct
realizations of $p$ if and only if the relation $I_p$ is
antisymmetric. Since $I_p$ is reflexive, the definition of
$\nu(p)$ and Propositions 1.4, 4.3 imply

\medskip
{\bf Corollary 4.4.} {\em Let $p(x)$ be a complete type of a
theory $T$ having a model $\mathcal{M}_p$, $\nu(p)$ be a regular
labelling function. The following conditions are equivalent:

$(1)$ the relation $I_p$ {\rm (}on the set of realizations of $p$
in any model $\mathcal{M}\models T${\rm )} is a partial order;

$(2)$ the structure $\mathfrak P_{\nu(p)}$ is an almost
deterministic monoid and $\rho_{\nu(p)}\subseteq U^{\leq 0}$.

This partial order $I_p$ is identical if and only if
$\rho_{\nu(p)}=\{0\}$. If $I_p$ is not identical, it has infinite
chains.}

\medskip
{\bf Definition} \cite{SuLP, Su041, Su08, BSV}. A countable model
$\mathcal{M}$ of theory $T$ is {\em limit} (accordingly {\em limit
over a type $p\in S(T)$}) if $\mathcal{M}$ is not prime over
tuples and
$\mathcal{M}=\bigcup\limits_{n\in\omega}\mathcal{M}({\bar{a}_n}),$
where $(\mathcal{M}({\bar{a}_n}))_{n\in\omega}$ is an elementary
chain of prime models over tuples $\bar{a}_n$ (and
$\mathcal{M}\models p(\bar{a}_n)$), $n\in\omega$.

\medskip
A characterization for the (non)symmetry of a relation $I_q$ for
the class of small theories is obtained in \cite{BSV}:

\medskip
{\bf Theorem 4.5.} {\em Let $q(\bar{x})$ be a complete type of a
small theory $T$. The following conditions are equivalent:

$(1)$ there exists a limit model over $q$;

$(2)$ the relation $I_q$ of isolation on a set of realizations of
$q$ in a {\rm (}any{\rm )} model $\mathcal{M}\models T$ realizing
$q$ is non-symmetric;

$(3)$ in some {\rm (}any{\rm )} model $\mathcal{M}\models T$
realizing $q$, there exist realizations $\bar{a}$ and $\bar{b}$ of
$q$ such that the type ${\rm tp}(\bar{b}/\bar{a})$ is principal
and $\bar{b}$ does not semi-isolate $\bar{a}$ and, in particular,
${\rm SI}_q$ is non-symmetric on $\mathcal{M}$.}

\medskip
Proposition 4.3 and Theorem 4.5 imply

\medskip
{\bf Corollary 4.6.} {\em Let $p(x)$ be a complete type of a small
theory $T$, $\nu(p)$ be a regular labelling function. ~The
following conditions are equivalent:

$(1)$ $I_p$ {\rm (}on the set of realizations of $p$ in any model
$\mathcal{M}\models T${\rm )} is an equivalence relation;

$(2)$ the structure $\mathfrak P_{\nu(p)}$ is an almost
deterministic monoid and there are no limit models over $p$;

$(3)$ the structure $\mathfrak P_{\nu(p)}$ is an almost
deterministic monoid and consists of non-negative labels.}

\medskip
In Corollary 4.6, the equivalence of (1) and (3) is implied by the
existence of $\mathcal{M}_p$ without the assumption of smallness
of $T$.

\medskip
{\bf Definition.} An element $u\in \rho_{\nu(p)}$ is called {\em
{\rm (}almost{\em )}
deterministic}\index{Element!deterministic}\index{Element!almost
deterministic} if for any/some realization $a$ of $p$ the formula
$\theta_u(a,y)$ has unique solution (has finitely many solutions).

\medskip
Note that there are no negative almost deterministic elements $u$
for a theory $T$ having an atomic model and finitely many
non-principal $1$-types in $S(T)$.\footnote{The following
arguments, in fact, repeat the remark after the proof of
Proposition 1.4.2 in \cite{SuLP}.} Indeed, otherwise the presence
of a negative element $u$ implies that the type $p(x)$ is
non-principal and the relation ${\rm SI}_p$ is not symmetric that
witnessed by the formula $\theta_u(x,y)$. Since for $\models p(a)$
the isolating formula $\theta_u(a,y)$ has some
$k\in\omega\setminus\{0\}$ solutions, there exists a formula
$\varphi(x)\in p(x)$ such that for any realization $b$ of
$\varphi(x)$ there are exactly $k$ solutions of the formula
$\theta_u(b,y)$. Moreover, since there are finitely many
non-principal $1$-types, there exist an element $c$, realizing a
principal type, and an element $d$ such that
$\models\varphi(c)\wedge\theta_u(c,d)\wedge\theta_u(a,d)$. It
means that the non-principal type $p$ is realized in an atomic
model that is impossible

At the same time, Example 1.4.3 in \cite{SuLP} illustrates that
there are theories $T$ with even deterministic negative elements
$u$, where there are infinitely many non-principal $1$-types in
$S(T)$.

\medskip
{\bf Proposition 4.7.} {\em If elements $u$ and $v$ are {\rm
(}almost{\rm )} deterministic then any element $v'$ in $u\cdot v$
is {\rm (}almost{\rm )} deterministic.}

\medskip
{\em Proof.} Consider formulas $\theta_u(a,y)$, $\theta_v(a,y)$,
and $\theta_{u,v}(a,y)$, where $\models p(a)$. If $u$ and $v$ are
deterministic then all these formulas have unique solutions, so
the element $v'\in u\cdot v$ is unique, and the formulas
$\theta_{u,v}(a,y)$ and $\theta_{v'}(a,y)$ are equivalent.

If $u$ and $v$ are almost deterministic then the formulas
$\theta_u(a,y)$, $\theta_v(a,y)$, and $\theta_{u,v}(a,y)$ have
finitely many solutions. It implies that the set $u\cdot v$ is
finite and there are finitely many solutions for the formulas
$\theta_{v'}(a,y)$, $v'\in u\cdot v$.~$\Box$

\medskip
Proposition 4.7 immediately implies

\medskip
{\bf Corollary 4.8.} {\em For any groupoid $\mathfrak P_{\nu(p)}$
its restriction $\mathfrak P_{\nu(p),d}$\index{$\mathfrak
P_{\nu(p),d}$} {\rm (}respectively $\mathfrak P_{\nu(p),{\rm
ad}}${\rm )}\index{$\mathfrak P_{\nu(p),{\rm ad}}$} to the set of
{\rm (}almost{\rm )} deterministic elements is a monoid too.}

\medskip
The following proposition resents a characterization for the
determinacy of non-negative elements in $\mathfrak P_{\nu(p)}$
assuming the presence of the model $\mathcal{M}_p$.

\medskip
{\bf Proposition 4.9.} {\em If the model $\mathcal{M}_p$ exists
then an element $u\geq 0$ in $\mathfrak P_{\nu(p)}$ is
deterministic if and only if $u^{-1}\cdot u=\{0\}$.}

\medskip
{\em Proof.} Let an element $u$ be deterministic, i.~e.,
$\theta_u(a,\mathcal{M}_p)=\{b\}$ for some realizations $a$ and
$b$ of $p$ in $\mathcal{M}_p$. Then
$\theta_{u^{-1},u}(b,\mathcal{M}_p)=\{b\}$, i.~e., $u^{-1}\cdot
u=\{0\}$.

We assume now that $u^{-1}\cdot u=\{0\}$ and prove that the
formula $\theta_u(a,y)$, where $\models p(a)$, has the unique
solution. Assume on the contrary that there are at least two
solutions $b_1$ and $b_2$. Then we have
$\models\theta_{u^{-1}}(b_1,a)\wedge\theta_u(a,b_2)$. Since $0\in
u^{-1}\cdot u$, $\theta_0(b_1,y)=(b_1\approx y)$, and
$\theta_0(b_1,y)\vdash\theta_{u^{-1},u}(b_1,y)$ then the
consistency of the formula
$\theta_{u^{-1},u}(b_1,y)\wedge\neg\theta_0(b_1,y)$ and the
existence of $\mathcal{M}_p$ imply that there is an isolating
formula $\theta_v(b_1,y)$, $v\ne 0$, such that
$\theta_v(b_1,y)\vdash\theta_{u^{-1},u}(b_1,y)$. It contradicts
the condition $u^{-1}\cdot u=\{0\}$.~$\Box$

\medskip
Unlike the determinacy there are no similar characterizations for
the almost determinacy.

\medskip
{\bf Example 4.2.} If $\Gamma=\langle M;R\rangle$ is an acyclic
undirected graph consisting of vertices of fixed degree $\upsilon$
then for the unique 1-type $p(x)\in S({\rm Th}(\Gamma))$, for the
principal formulas $\theta_n(x,y)$, where
$\models\theta_n(a,b)\Leftrightarrow\rho(a,b)=n$, $n\in\omega$,
and for the monoid $\mathfrak P_{\nu(p)}$ over the alphabet
$\omega$ we have $m\cdot n=\{m+n,|m-n|\}$. In particular,
$n=n^{-1}$ and $n\cdot n=\{0,2n\}$. At the same time the monoid
$\mathfrak P_{\nu(p)}$ does not depend on
$\upsilon\in\omega\cup\{\infty\}$.

\medskip
{\bf Proposition 4.10.} {\em If $\mathfrak P_{\nu(p)}$ is a
deterministic monoid then the structure $\mathfrak P'_{\nu(p)}$ is
a group if and only if $\rho_{\nu(p)}$ consists of non-negative
elements.}

\medskip
{\em Proof.} At first we observe that, by definition, if
$u\in\rho_{\nu(p)}$ is negative then there are no labels $v$ such
that $u\odot v=0$. Hence, if ${\rm PFN}(p)\ne\varnothing$ then
$\mathfrak P'_{\nu(p)}$ is not a group.

Now we assume that $\rho_{\nu(p)}\cap U^-=\varnothing$ and prove
that the structure $\mathfrak P'_{\nu(p)}$ is a group. Indeed, if
${\rm PFN}(p)=\varnothing$ then for any element
$u\in\rho_{\nu(p)}$ there is the (unique) inverse element
$v=u^{-1}$ such that $0\in u\cdot v$. As the monoid $\mathfrak
P_{\nu(p)}$ is deterministic we obtain $u\odot v=0$.~$\Box$

\medskip
{\bf Corollary 4.11.} {\em If the model $\mathcal{M}_p$ exists,
the monoid $\mathfrak P_{\nu(p)}$ is deterministic, and $\mathfrak
P'_{\nu(p)}$ is a group, then all elements in $\mathfrak
P'_{\nu(p)}$ are deterministic.}

\medskip
{\em Proof.} Since by Proposition 4.10 the set $\rho_{\nu(p)}$
consists of non-negative elements then, as the monoid $\mathfrak
P_{\nu(p)}$ is deterministic, by Proposition 4.9 each element in
$\mathfrak P'_{\nu(p)}$ is deterministic.~$\Box$

\medskip
{\bf Proposition 4.12.} {\em If the model $\mathcal{M}_p$ exists
then the set $\rho_{\nu(p),d}^{\geq 0}$ of all non-negative
deterministic elements $u$ in $\rho_{\nu(p)}$, for which elements
$u^{-1}$ are also deterministic, forms a deterministic submonoid
$\mathfrak G_{\nu(p)}$\index{$\mathfrak G_{\nu(p)}$} of the monoid
$\mathfrak P_{\nu(p),d}$, consisting of deterministic elements of
$\mathfrak P_{\nu(p)}$, and such that $(\mathfrak G_{\nu(p)})'$ is
a group.}

\medskip
{\em Proof.} Since for any $u\in \rho_{\nu(p),d}^{\geq 0}$ the
element $u^{-1}$ satisfying $u\cdot u^{-1}=u^{-1}\cdot u=\{0\}$
belongs to $\rho_{\nu(p),d}^{\geq 0}$ it suffices to observe that
if $u,v\in \rho_{\nu(p),d}^{\geq 0}$ then $u\cdot v$ contains a
unique element $v'$ and this element is deterministic by
Proposition 4.7.~$\Box$

\medskip
In Figure 1, a Hasse diagram is presented illustrating the links
of the structure $\mathfrak P_{\nu(p)}$ with structures above,
being restrictions of $\mathfrak P_{\nu(p)}$ to subalphabets of
$U$. Here the superscripts $\cdot^{\leq 0}$ and $\cdot^{\geq 0}$
point out on restrictions of $\mathfrak P_{\nu(p)}$ to the sets of
non-positive and non-negative elements respectively, the
subscripts $\cdot_d$ and $\cdot_{\rm ad}$ indicate the sets of
deterministic and almost deterministic elements. By Propositions
3.1 and 4.2, just $\mathfrak P_{\nu(p)}$ and $\mathfrak
P_{\nu(p)}^{\leq 0}$ may not be monoids.

\begin{figure}[t]
\begin{center}
\unitlength 2cm
\begin{picture}(8,4.5)(0.65,-1.3)
{\small

\put(3,0){\makebox(0,0)[cc]{$\bullet$}}
\put(3,1){\makebox(0,0)[cc]{$\bullet$}}
\put(3,2){\makebox(0,0)[cc]{$\bullet$}}
\put(4,-1){\makebox(0,0)[cc]{$\bullet$}}
\put(4,1){\makebox(0,0)[cc]{$\bullet$}}
\put(4,2){\makebox(0,0)[cc]{$\bullet$}}
\put(4,3){\makebox(0,0)[cc]{$\bullet$}}
\put(5,-0.5){\makebox(0,0)[cc]{$\bullet$}}
\put(5,0){\makebox(0,0)[cc]{$\bullet$}}
\put(5,1){\makebox(0,0)[cc]{$\bullet$}}
\put(5,2){\makebox(0,0)[cc]{$\bullet$}} \put(3,1){\line(0,1){1}}
\put(3,1){\line(0,-1){1}} \put(3,2){\line(1,1){1}}
\put(5,1){\line(0,-1){0.5}} \put(5,1){\line(0,1){1}}
\put(4,0){\line(0,1){1}} \put(4,2){\line(0,1){1}}
\put(5,1){\line(0,-1){1}} \put(4,1){\line(0,-1){1}}
\put(4,1){\line(0,1){1}} \put(4,0){\line(0,-1){1}}
\put(4,1){\line(-1,-1){1}} \put(4,1){\line(1,-1){1}}
\put(4,3){\line(0,-1){1}} \put(4,2){\line(-1,-1){1}}
\put(4,3){\line(1,-1){1}} \put(4,2){\line(1,-1){1}}
\put(5.1,-0.5){\makebox(0,0)[cl]{$\mathfrak G_{\nu(p)}$}}

\put(3,0){\line(1,-1){1}}

\put(5,-0.5){\line(0,1){1}}

\put(5,-0.5){\line(-2,-1){1}}

\put(3.13,0){\makebox(0,0)[cl]{$\mathfrak P^{\leq 0}_{\nu(p),d}$}}
\put(4.1,-1.1){\makebox(0,0)[cl]{$\{0\}$}}
\put(5.1,0){\makebox(0,0)[cl]{$\mathfrak P^{\geq 0}_{\nu(p),d}$}}
\put(3.13,1){\makebox(0,0)[cl]{$\mathfrak P^{\leq 0}_{\nu(p),{\rm
ad}}$}} \put(4.1,1){\makebox(0,0)[cl]{$\mathfrak P_{\nu(p),d}$}}
\put(5.1,1){\makebox(0,0)[cl]{$\mathfrak P^{\geq 0}_{\nu(p),{\rm
ad}}$}} \put(3.13,2){\makebox(0,0)[cl]{$\mathfrak P^{\leq
0}_{\nu(p)}$}} \put(4.1,2){\makebox(0,0)[cl]{$\mathfrak
P_{\nu(p),{\rm ad}}$}} \put(5.1,2){\makebox(0,0)[cl]{$\mathfrak
P^{\geq 0}_{\nu(p)}$}} \put(4.1,3.1){\makebox(0,0)[cl]{$\mathfrak
P_{\nu(p)}$}}}

\end{picture}

\vspace{1cm} Fig. 1
\end{center}
\end{figure}

\medskip
\section{Graph compositions \\ and monoid compositions}
\medskip

Recall \cite{SO1} that the \emph{composition}\index{Composition!of
graphs} $\Gamma_1[\Gamma_2]$\index{$\Gamma_1[\Gamma_2]$} {\em of
graphs} $\Gamma_1=\langle X_1; R_1\rangle$ and $\Gamma_2=\langle
X_2; R_2\rangle$ is the graph $\langle X_1\times X_2;R\rangle$,
where $((a_1,b_1),(a_2,b_2))\in R$ if and only if some of the
following conditions is met:

1) $(a_1,a_2)\in R_1$;

2) $a_1=a_2$ and $(b_1,b_2)\in R_2$.

Similarly we define the notion of monoid composition.

Let $\mathcal{S}_1$ and $\mathcal{S}_2$ be monoids, for which $0$
is the unit, $S_1\subseteq U^{\leq 0}$, and $S_2\subseteq U^{\geq
0}$. The {\em composition}\index{Composition!of monoids} or the
{\em successively-annihilating
band}\index{Band!successively-annihilating}\footnote{see
\cite{Lyap, OA2}.}
$\mathcal{S}_1[\mathcal{S}_2]$\index{$\mathcal{S}_1[\mathcal{S}_2]$}
{\em of monoids} $\mathcal{S}_1$ and $\mathcal{S}_2$ is the
algebra $\langle S_1\cup S_2;\,\odot\rangle$, where $\langle
S_1\cup S_2;\,\odot\rangle\upharpoonright S_i=\mathcal{S}_i$,
$i=1,2$, and $u\odot v=v\odot u=u$ for $u<0$ and $v>0$.

\medskip
{\bf Proposition 5.1} \cite{Lyap}. {\em Any
successively-annihilating band $\mathcal{S}_1[\mathcal{S}_2]$ is a
monoid.}

\medskip
{\em Proof.} We fix a  successively-annihilating band
$\mathcal{S}_1[\mathcal{S}_2]$. Since $\mathcal{S}_1$ and
$\mathcal{S}_2$ are monoids it suffices to check the
associativity: $(u_1\cdot u_2)\cdot u_3=u_1\cdot (u_2\cdot u_3)$
for any three elements $u_1,u_2,u_3$, where exactly two of them
belong to $U^-$ or to $U^+$.

We check this property analyzing six cases:

(1) if $u_1\in U^-$ and $u_2,u_3\in U^+$ then $$(u_1\odot
u_2)\odot u_3=u_1\odot u_3=u_1=u_1\odot (u_2\odot u_3);$$

(2) if  $u_1, u_3\in U^+$ and $u_2\in U^-$ then $$(u_1\odot
u_2)\odot u_3=u_2\odot u_3=u_2=u_1\odot u_2=u_1\odot (u_2\odot
u_3);$$

(3) if $u_1, u_2\in U^+$ and $u_3\in U^-$ then $$(u_1\odot
u_2)\odot u_3=(u_1\odot u_2)\odot u_3=u_3=u_1\odot u_3=u_1\odot
(u_2\odot u_3);$$

(4) if $u_1,u_2\in U^-$ and $u_3\in U^+$ then $$(u_1\odot
u_2)\odot u_3=u_1\odot u_2=u_1\odot(u_2\odot u_3)=u_1\odot
(u_2\odot u_3);$$

(5) if $u_1,u_3\in U^-$ and $u_2\in U^+$ then $$(u_1\odot
u_2)\odot u_3=u_1\odot u_3=u_1\odot(u_2\odot u_3)=u_1\odot
(u_2\odot u_3);$$

(6) if $u_1\in U^+$ and $u_2,u_3\in U^-$ then $$(u_1\odot
u_2)\odot u_3=u_2\odot u_3=u_1\odot(u_2\odot u_3)=u_1\odot
(u_2\odot u_3).\,\,\Box$$

\medskip
{\bf Theorem 5.2.} {\em For any group $\langle G;\,\ast\rangle$,
where the universe consists of non-negative elements and $0$
denotes the group unit, and for the monoid
$\langle\{-1,0\};\,+\rangle$ with the zero element $0$ and the
idempotent element $-1$ there is a theory $T$ with a type $p\in
S(T)$ and a regular labelling function $\nu(p)$ such that the
monoid $\mathfrak P'_{\nu(p)}$ coincides with the monoid
$\langle\{-1,0\};\,+\rangle[\langle G;\,\ast\rangle]$.}

\medskip
{\em Proof.} We construct a structure $\mathcal{M}$ such that its
theory $T={\rm Th}(\mathcal{M})$ has a type $p(x)\in S(T)$ and a
regular labelling function $\nu(p)$ with $\mathfrak
P'_{\nu(p)}=\langle\{-1,0\};\,+\rangle[\langle G;\,\ast\rangle]$.
For this aim we consider the Ehrenfeucht's example $\langle\mathbb
Q;<,c_k\rangle_{k\in\omega}$, $c_k<c_{k+1}$, $k\in\omega$, such
that each element $a$ is replaced by a $<$-antichain consisting of
$|G|$ elements and forming a free $1$-generated polygon over the
group $\langle G;\,\ast\rangle$ isomorphic to the structure
$\mathcal{G}=\langle G; Q_g\rangle_{g\in G}$, where
$Q_g=\{(a,b)\in G^2\mid a\ast g=b\}$, $g\in G$. Here we replace
each constant $c_k$ by a unary predicate $R_k$ consisting of
elements of a copy of $\mathcal{G}$. Thus we form the composition
$\langle\mathbb Q;<\rangle[\mathcal{G}]$ of graphs expanded by
relations $R_k$, $k\in\omega$, $(x<y)$,
$\neg(x<y)\wedge\neg(y<x)$, $Q_g$, $g\in G$. The unique
non-principal $1$-type $p(x)$ is isolated by set of formulas
$\exists y(R_k(y)\wedge(y<x))$, $k\in\omega$. For any realization
$a$ of $p$ the list of pairwise non-equivalent isolating formulas
$\varphi(a,y)$ with $\varphi(a,y)\vdash p(y)$ is exhausted by the
formulas $(a<y)$ and $Q_g(a,y)$, $g\in G$. We define a regular
labelling function $\nu(p)$ such that the formula $(a<y)$ has the
label $-1$ and the formulas $Q_g(a,y)$ have non-negative labels
$g$. Since $<\circ <\,\,=\,\,<$, $<\,\circ\,\,
Q_g=Q_g\,\,\circ\,<\,\,=\,\,<$, $g\in G$, and the links between
elements of $\rho_{\nu(p)}^{\geq 0}$ are defined by the group
$\langle G;\,\ast\rangle$, the monoid $\mathfrak P'_{\nu(p)}$
coincides with the monoid $\langle\{-1,0\};\,+\rangle[\langle
G;\,\ast\rangle]$.~$\Box$

\medskip
{\bf Theorem 5.3.} {\em For any group $\langle G;\,\ast\rangle$
consisting of non-negative elements with the unit element $0$ and
for the monoid $\langle\omega^\ast;\,+\rangle$ of non-positive
integers there exists a theory $T$ with a type $p\in S(T)$ and a
regular labelling function $\nu(p)$ such that the monoid
$\mathfrak P'_{\nu(p)}$ coincides with the monoid
$\langle\omega^\ast;\,+\rangle[\langle G;\,\ast\rangle]$.}

\medskip
{\em Proof.} We construct a structure $\mathcal{M}$ such that its
theory $T={\rm Th}(\mathcal{M})$ has a type $p(x)\in S(T)$ and a
regular labelling function $\nu(p)$ with $\mathfrak
P'_{\nu(p)}=\langle\omega^\ast;\,+\rangle[\langle
G;\,\ast\rangle]$.

The language of $\mathcal{M}$ consists of unary predicate symbols
${\rm Col}_n$, $n\in\omega$ (forming a coloring of the set $M$) of
binary predicate symbol $Q$, and of binary predicate symbols
$Q_g$, $g\in G$.

We consider a connected acyclic directed graph $\Gamma=\langle
M_0;Q\rangle$, where each element has infinitely many images and
infinitely many preimages, i.~e., $\Gamma$ forms a free directed
pseudoplane \cite{SuLP, Pi89_2, Su901}.

We define an $1$-inessential $Q$-ordered coloring (see
\cite{SuLP}) ${\rm Col }\mbox{\rm : }M_0\to\omega\cup\{\infty\}$
of $\Gamma$ producing unary predicates ${\rm Col}_n=\{a\in M_0\mid
{\rm Col}(a)=n\}$, $n\in\omega$.

For the graph  $\Gamma$ we define, by induction, relations
$Q^n$\index{$Q^n$}, $n\in{\bf Z}$: $Q^0\rightleftharpoons{\rm
id}_{M_0}$, $Q^{n+1}\rightleftharpoons Q^n\circ Q$,
$Q^{-n}\rightleftharpoons (Q^n)^{-1}$, $n\in\omega$.

Note that for the (unique) non-principal type $p(x)$, isolated by
the set $\{\neg{\rm Col }_m(x)\mid m<\omega\}$ of formulas, and
for any realizations $a$ and $b$ of $p$, the pair $(a,b)$ is a
principal arc if and only if $\models Q^n(a,b)$ for some
$n\in\omega$.

We assume that the formula $Q^n(x,y)$ has the label $-n\in U^{\leq
0}$, $n\in\omega$. Since for any $m,n\in\omega$ the formula
$\exists z(Q^m(x,z)\wedge Q^n(z,y))$ is equivalent to the formula
$Q^{m+n}(x,y)$, then for the $Q$-structure on a set of
realizations of $p$ the structure $\mathfrak P'_{\nu(p)}$
coincides with $\langle\omega^\ast;\,+\rangle$.

Now we consider the group $\langle G;\,\ast\rangle$ and define, on
the set $G$, binary predicates $Q_g$, $g\in G$, by the rule:
$$
Q_g=\{(a,b)\in G^2\mid a\ast g=b\}.
$$
As in the proof of Theorem 5.2 the structure $\mathcal{G}=\langle
G; Q_g\rangle_{g\in G}$ forms a free $1$-generated polygon over
the group $\langle G;\,\ast\rangle$.

We define a model of required theory $T$ as the composition
$\Gamma[\mathcal{G}]$ of graphs with colored vertices and arcs
such that each vertex $a$ of $\Gamma$ is replaced by a copy of
structure $\mathcal{G}$, for which all elements have the color
${\rm Col}(a)$. The relations $Q_g$, for $\Gamma[\mathcal{G}]$,
are composed as the unions of correspondent relations in the
copies of $\mathcal{G}$, and the relation $Q$, in
$\Gamma[\mathcal{G}]$, consists of all pairs $(a',b')$, where
$a'\in C_a$, $b'\in C_b$, $(a,b)\in Q$ in $\Gamma$, and $C_a$,
$C_b$ are copies of $\mathcal{G}$ replacing vertices $a,b\in M_0$.
The composition preserves the uniqueness of the non-principal type
$p(x)$.

It remains to note that for any realization $a$ of $p$ the list of
pairwise non-equivalent isolating formulas $\varphi(a,y)$ with
$\varphi(a,y)\vdash p(y)$ is exhausted by the formulas $Q^n(a,y)$,
$n\in\omega$, $Q_g(a,y)$, $g\in G$, we have $\mathfrak
P'_{\nu(p)}\upharpoonright\omega^\ast=\langle\omega^\ast;\,+\rangle$,
$\mathfrak P'_{\nu(p)}\upharpoonright G=\langle G;\,\ast\rangle$
and $Q^n\circ Q_g=Q_g\circ Q^n=Q^n$ for $n>0$, $g\in G$.~$\Box$

\medskip
\section{$I$-groupoids}
\medskip

In this section, we collect basic structural properties of
$I_{\nu(p)}$-groupoids and prove that any groupoid $\mathfrak P$
satisfying that list of properties coincides with some
$I_{\nu(p)}$-monoid $\mathfrak P_{\nu(p)}$.

Let $U=U^-\,\dot{\cup}\,\{0\}\,\dot{\cup}\,U^+$ be an alphabet
consisting of a set $U^-$\index{$U^-$} of {\em negative
elements}\index{Element!negative}, a set $U^+$\index{$U^+$} of
{\em positive elements}\index{Element!positive} and a zero $0$. As
above we write $u<0$ for any element $u\in U^-$, $u>0$ for any
element $u\in U^+$, and $u\cdot v$ instead of $\{u\}\cdot\{v\}$
considering an operation $\cdot$ on the set
$\mathcal{P}(U)\setminus\{\varnothing\}$.

A groupoid $\mathfrak
P=\langle\mathcal{P}(U)\setminus\{\varnothing\};\,\cdot\rangle$ is
called an {\em $I$-groupoid}\index{$I$-groupoid} if it satisfies
the following conditions:

\medskip
${\small\bullet}$ the set $\{0\}$ is the unit of the groupoid
$\mathfrak P$;

\medskip
${\small\bullet}$ the operation $\cdot$ of the groupoid $\mathfrak
P$ is generated by the function $\cdot$ on elements in $U$ such
that each elements $u,v\in U$ define a nonempty set $(u\cdot
v)\subseteq U$: for any sets
$X,Y\in\mathcal{P}(U)\setminus\{\varnothing\}$ the following
equality holds:
$$
X\cdot Y=\bigcup\{x\cdot y\mid x\in X,y\in Y\};
$$

\medskip
${\small\bullet}$ if $u<0$ then the sets $u\cdot v$ and $v\cdot u$
consist of negative elements for any $v\in U$;

\medskip
${\small\bullet}$ if $u>0$ and $v>0$ then the set $u\cdot v$
consists of non-negative elements;

\medskip
${\small\bullet}$ for any $u>0$ there is the unique {\em
inverse}\index{Element!inverse} element $u^{-1}>0$ such that $0\in
(u\cdot u^{-1})\cap(u^{-1}\cdot u)$;

\medskip
${\small\bullet}$ is a positive element $u$ belongs to a set
$v_1\cdot v_2$ then $u^{-1}$ belongs to $v_2^{-1}\cdot v_1^{-1}$;

\medskip
${\small\bullet}$ for any elements $u_1,u_2,u_3\in U$ the
following inclusion holds: $$(u_1\cdot u_2)\cdot u_3\supseteq
u_1\cdot(u_2\cdot u_3),$$ and the strict inclusion
$$(u_1\cdot u_2)\cdot u_3\supset u_1\cdot(u_2\cdot u_3)$$ may be
satisfied only for $u_1<0$ and $|u_2\cdot u_3|\geq\omega$;

\medskip
${\small\bullet}$ the groupoid $\mathfrak P$ contains the {\em
deterministic}\index{Subgroupoid!deterministic} subgroupoid
$\mathfrak P^{\geq 0}_d$\index{$\mathfrak P^{\geq 0}_d$} (being a
monoid) with the universe $\mathcal{P}(U^{\geq
0}_d)\setminus\{\varnothing\}$, where $$U^{\geq 0}_d=\{u\in
U^{\geq 0}\mid u^{-1}\cdot u=\{0\}\};$$ any set $u\cdot v$ is a
singleton, where $u,v\in U^{\geq 0}_d$.

\medskip
By the definition each $I$-groupoid $\mathfrak P$ contains
$I$-subgroupoids $\mathfrak P^{\leq 0}$\index{$\mathfrak P^{\leq
0}$} and $\mathfrak P^{\geq 0}$\index{$\mathfrak P^{\geq 0}$} with
the universes $\mathcal{P}(U^-\cup\{0\})\setminus\{\varnothing\}$
and $\mathcal{P}(U^+\cup\{0\})\setminus\{\varnothing\}$
respectively. The structure $\mathfrak P^{\geq 0}$ is a monoid.

\medskip
{\bf Theorem 6.1.} {\em For any {\rm (}at most countable{\rm )}
$I$-groupoid $\mathfrak P$ there is a {\rm (}small{\rm )} theory
$T$ with a type $p(x)\in S(T)$ and a regular labelling function
$\nu(p)$ such that $\mathfrak P_{\nu(p)}=\mathfrak P$.}

\medskip
{\em Proof.} We fix an $I$-monoid $\mathfrak
P=\langle\mathcal{P}(U)\setminus\{\varnothing\};\,\cdot\rangle$.
The construction of a required theory will be fulfilled in
accordance with a construction of a generic structure
$\mathcal{M}$ of language $\Sigma=\{{\rm Col}^{(1)}_n\mid
n\in\omega\}\cup\{Q^{(2)}_u\mid u\in U\}$ \cite[Chapter 2]{SuLP}
with pairwise disjoint predicates $Q_u$, with an ordered coloring
${\rm Col}\mbox{\rm : }M\to\omega\cup\{\infty\}$ with respect to
each formula $Q_u(x,y)$, where $u<0$, and with a unique
non-principal $1$-type $p(x)$ (isolated by the set $\{\neg{\rm
Col}_n(x)\mid n\in\omega\}$ of formulas). W.l.o.g. we assume that
$|U|\leq\omega$ (for $|U|>\omega$, the construction differs by
cardinalities of diagrams describing links for elements of finite
sets and by cardinalities of sets of diagrams forming generic
models).

Consider a generic class $({\bf T}_0;\leqslant)$ consisting of all
possible diagrams $\Phi(A)$ over finite sets $A$ such that each
$\Phi(A)$ contains a maximal consistent set of quantifier-free
formulas $\varphi(\bar{a})$, $\bar{a}\in A$, united with a set of
formulas $Q^\delta_{uv}(a,b)$, $a,b\in A$, $\delta\in\{0,1\}$,
$Q_{uv}(x,y)=\exists z(Q_{u}(x,z)\wedge Q_{v}(z,y))$, $u,v\in U$,
and $\Phi(A)$ includes formulas with parameters in $A$, without
free variables, and describing the following properties:

\medskip
(1) for any $u\in U$ any element in $A$ is a image and a preimage
of some elements by the relation $Q_u$;

(2) the relation $Q_0$ on the set $A$ is identical;

(3) if $a\in A$ then all $Q_u$-images of $a$ have colors $\geq
{\rm Col}(a)$ and all $Q_u$-preimages of $a$ have the colors $\leq
{\rm Col}(a)$;

(4) if $u>0$, $a\in A$, and $Q_u(a,b)\in\Phi(A)$ then
$Q_{u^{-1}}(b,a)\in\Phi(A)$ and ${\rm Col}(b)={\rm Col}(a)$;

(5) if $v\in u_1\cdot u_2$ and $Q_v(a,b)\in\Phi(A)$ then
$Q_{u_1u_2}(a,b)\in\Phi(A)$;

(6) for any $u\ne 0$ some diagram $\Psi(B)\supseteq\Phi(A)$ in
${\bf T}_0$ defines a graph $\langle B;Q_u\rangle$ with a cycle if
and only if $0\in \underbrace{u\cdot\ldots\cdot
u}_{n\mbox{\scriptsize \ times}}$ for some $n>0$;

(7) if $u\in U^{\geq 0}_d$ then each element $a\in A$ has a unique
$Q_u$-image; the following inductive condition describes the least
set $U^{\geq 0}_{\rm ad}\supseteq U^{\geq 0}_d$ of non-negative
elements $u\in U$ for which the sets of $Q_u$-images and of
$Q_u$-preimages of $a$ are finite: if $(u\cdot
u^{-1})\cup(u^{-1}\cdot u)$ consists of finitely many elements
belonging to $U^{\geq 0}_{\rm ad}$ then $u,u^{-1}\in U^{\geq
0}_{\rm ad}$; if $u^{-1}\cdot u$ consists of finitely many
elements belonging to $U^{\geq 0}_{\rm ad}$ then each element $a$
has finitely many $Q_u$-images; if $u\cdot u^{-1}$ consists of
finitely many elements belonging to $U^{\geq 0}_{\rm ad}$ then
each element $a$ has finitely many $Q_u$-preimages; for other
elements $u$ the numbers of $Q_u$-images and of $Q_u$-preimages
for elements $a\in A$ is unbounded;

(8) if $u_1,u_2\in U$ and the set $u_1\cdot u_2$ is (in)finite
then for any element $a\in A$ the set of $Q_{u_1u_2}$-images of
$a$ is represented as a union of sets of $Q_v$-images for all
elements $v\in u_1\cdot u_2$ (and some set of elements that are
not $Q_u$-images of $a$ on any of the relations $Q_u$);

(9) for any element $v\in((u_1\cdot u_2)\cdot u_3)\setminus(
u_1\cdot(u_2\cdot u_3))$ there is a description forming Example
3.1.

\medskip
If $\Phi(A),\Psi(B)$ are diagrams in ${\bf T}_0$ and
$\Phi(A)\subseteq\Psi(B)$, we suppose, by the definition, that
$\Phi(A)$ is a {\em strong subdiagram}\index{Subdiagram!strong} of
$\Psi(B)$ (i.~e.,
$\Phi(A)\leqslant\Psi(B)$)\index{$\Phi(A)\leqslant\Psi(B)$} if
$A$, with each element $a$ in $A$, contains all its $Q_u$-images
in $B$, where $u^{-1}\cdot u$ consists of finitely many labels
belonging to $U^{\geq 0}_{\rm ad}$.

For the checking that $({\bf T}_0;\leqslant)$ is a self-sufficient
generic class, it suffices to observe that for any diagrams
$\Phi(A),\Psi(B),$ ${\rm X}(C)\in{\bf T}_0$ with
$\Phi(A)\leqslant\Psi(B)$, $\Phi(A)\leqslant{\rm X}(C)$, and
$A=B\cap C$ there is a diagram $\Theta(B\cup C)\in{\bf T}_0$ such
that $\Psi(B)\leqslant\Theta(B\cup C)$ and ${\rm
X}(C)\leqslant\Theta(B\cup C)$.

For the type $\Theta(B\cup C)$ we choose the set $\Psi(B)\cup{\rm
X}(C)$ extended by the following formulas for elements $b\in
B\setminus A$ and $c\in C\setminus A$:

\medskip
(a) $\theta_{u,v}(b,c)$, where $Q_u(b,a)\in\Psi(B)$ and
$Q_v(a,c)\in{\rm X}(C)$ for some $a\in A$;

(b) $\neg\theta_{u,v}(b,c)$, where $\neg Q_u(b,a)\in\Psi(B)$ or
$\neg Q_v(a,c)\in{\rm X}(C)$ for all $a\in A$;

(c) some formulas $\theta_{v'}(b,c)$, where $Q_u(b,a)\in\Psi(B)$
and $Q_v(a,c)\in{\rm X}(C)$ for some $a\in A$, $v'\in u\cdot v$,
and the set $u\cdot v$ is finite;

(d) formulas $\neg\theta_{v'}(b,c)$, $v'\in U$, if the previous
items do not imply a converse.

\medskip
We claim that, applying the generic construction, one obtains a
$({\bf T}_0;\leqslant$~$)$-generic saturated structure
$\mathcal{M}$ with the generic theory $T={\rm Th}(\mathcal{M})$,
the type $p(x)\in S(T)$, and the regular labelling function
$\nu(p)\mbox{\rm : }{\rm PF}(p)/{\rm PE}(p)\to U$ satisfying the
condition $\mathfrak P_{\nu(p)}=\mathfrak P$. By Proposition
1.2.13 in \cite{SuLP}, each formula $Q_u(x,y)$, $u<0$, witnesses
on the non-symmetry of the relation ${\rm SI}_p$, and each formula
$Q_u(x,y)$, $u>0$, links realizations of $p$ only with
realizations of the same type and, being a principal formula of
the structure on the set $p(M)$ of realizations of $p$, has the
inverse principal formula $Q_{u^{-1}}(x,y)$ on $p(M)$.

Now we argue to show that $\mathcal{M}$ is saturated. If $U^{\geq
0}_{\rm ad}$ is finite the saturation of $\mathcal{M}$ is implied
by \cite[Theorem 2.5.1]{SuLP} (see also \cite[Theorem 4.1]{Su073})
in view of the uniform $t$-amalgamation property that holds by the
formula definability of self-sufficient closure of any finite set.

Using the proof of the same theorem, we shall observe that
$\mathcal{M}$ is saturated for $|U^{\geq 0}_{\rm ad}|=\omega$. For
this aim we enumerate all predicates $Q_u$, $u\in U$: $Q_m$,
$m\in\omega$.

Let $\mathcal{M}'$ be an $\omega$-saturated model of  ${\rm
Th}(\mathcal{M})$, $\Phi(A)$ and $\Phi(A')=[\Phi(A)]^A_{A'}$ be
diagrams in ${\bf T}_0$ such that $\mathcal{M}\models\Phi(A)$ and
$\mathcal{M}'\models\Phi(A')$. If $\Psi(B')\in{\bf T}_0$,
$\Phi(A')\leqslant\Psi(B')$, and $\mathcal{M}'\models\Psi(B')$
then the construction of~$\mathcal{M}$ implies that there exists a
set $B\subset M$ extending $A$ and satisfying
$\mathcal{M}\models\Psi(B)$. It means that for a partial
isomorphism $f\mbox{\rm : }A\to A'$ between $\mathcal{M}$ and
$\mathcal{M}'$ there exists a partial isomorphism $g\mbox{\rm :
}B\to B'$ between these structures extending $f$.

Now, let $\Psi(B)\in{\bf T}_0$, $\Phi(A)\leqslant\Psi(B)$,
$\mathcal{M}\models\Psi(B)$, and $X$ and~$Y$ be disjoint sets of
variables, which are in bijective correspondence with sets $A$ and
$B\setminus A$. Assume that the formula $\varphi_n(X)$
($\psi_n(X,Y)$, respectively), $n\in\omega$, describes the
following:

(i) finite colors of elements of $A$ (of $B$);

(ii) negations of colors not exceeding $n$ for elements of $A$ (of
$B$) that are infinite in color;

(iii) the existence, colors of arcs, the existence and colors of
some arcs of pathes of length $2$ (including all possibilities for
colors $\leq n$ of intermediate arcs) connecting elements of $A$
(of $B$), and the colors $m\leq n$ of arcs outgoing from vertices
$a\in A$ ($a\in B$) for which $\exists y Q_m(a,y)\in\Phi(A)$
($\exists y Q_m(a,y)\in\Psi(B)$), $Q_m=Q_u$, $u\in U^{\geq 0}_{\rm
ad}$;

(iv) the non-existence of arcs of colors $\leq n$ and of pathes of
length $2$ (including all possibilities for colors $\leq n$ of
intermediate arcs) connecting elements of $A$ (of $B$), if these
elements are not linked by the pathes, as well as the absence of
colors $m\leq n$ for arcs outgoing from vertices $a\in A$ ($a\in
B$) for which $\neg\exists y Q_m(a,y)\in\Phi(A)$ ($\neg\exists y
Q_m(a,y)\in\Psi(B)$), $Q_m=Q_u$, $u\in U^{\geq 0}_{\rm ad}$.

By the construction of $\mathcal{M}$,
$$
\mathcal{M}\models\forall X\:(\varphi_n(X)\to\exists
Y\:\psi_n(X,Y)).
$$
Hence
$$
\mathcal{M}'\models\forall X\:(\varphi_n(X)\to\exists
Y\:\psi_n(X,Y)).
$$
This implies that the set $\{\psi_n(A',Y)\mid n\in\omega\}$ of
formulas is locally realizable in $\mathcal{M}'$; hence, it is
realizable in $\mathcal{M}'$ since $\mathcal{M}'$ is
$\omega$-saturated. Therefore there exist a set $B'\subset M'$
containing $A'$, and a partial isomorphism $g\mbox{\rm : }B\to B'$
extending the partial isomorphism $f$.

The possibility for extending any partial isomorphisms $f\mbox{\rm
: }A\to A'$ and the known back-and-forth method show that the
structure $\mathcal{M}$ with distinguished constants for the
elements in $A\subset M$ is isomorphic to a countable elementary
substructure of the structure $\mathcal{M}'$ with distinguished
constants for the elements in $A'$. Since the finite sets $A$ and
$A'$ connected by a partial isomorphism and preserving a type
$\Phi(X)$ are chosen arbitrarily, and $\mathcal{M}'$ is saturated,
we conclude that $\mathcal{M}$ realizes any type over a finite
set, $\mathcal{M}$ is saturated, and ${\rm Th}(\mathcal{M})$ is
small.

Note that for the $({\bf T}_0;\leqslant$~$)$-generic structure
$\mathcal{M}$, the possibility for extending any finite partial
isomorphisms preserving types $\Phi(X)$ in ${\bf T}_0$ implies
that if~$A,B\subset M$, $\mathcal{M}\models\Phi(A)$ and
$\mathcal{M}\models\Phi(B)$ then there is an automorphism of
$\mathcal{M}$ extending the initial partial isomorphism between
$A$ and $B$. Consequently, ${\rm tp}_\mathcal{M}(A)={\rm
tp}_\mathcal{M}(B)$. In particular, for any realization $a$ of $p$
and for any $u\in U$ the formula $Q_u(a,y)$ is isolating and these
formulas exhaust the list of all pairwise non-equivalent isolating
formulas $\varphi(a,y)$ for which $\varphi(a,y)\vdash
p(y)$.~$\Box$

\medskip
{\bf Remark 6.2.} If an $I$-groupoid $\mathfrak P$ is constructed
by a set $U^{\geq 0}$ then by the construction above (restricting
the construction to a set of realizations of the type infinite in
color) there is a transitive theory $T$ with a (unique) type
$p(x)\in S(T)$ and a regular labelling function $\nu(p)$ such that
$\mathfrak P_{\nu(p)}=\mathfrak P$.

\medskip
\section{Groupoids of binary isolating formulas \\
on sets of realizations of types \\ of special theories}
\medskip

In this section, we present a specificity of groupoids
$\mathfrak{P}_{\nu(p)}$ for types $p$ of {\em
special}\index{Theory!special} theories used for the
classifications of countable models of Ehrenfeucht theories
\cite{SuLP, Su041, Su043}, of theories with finite Rudin--Keisler
preorders \cite{SuLP, Su072}, of small theories \cite{SuLP, Su08},
of $\omega$-stable theories with respect to numbers of limit
models over types \cite{Su10Irk}, as well as for the
investigations of graph links for limit models over types that
obtained by quotients of numerical sequences \cite{Shul, Shul2,
Shul3}.

\medskip
Let $\Gamma=\langle X,Q\rangle$ be a graph, and $a$ be a vertex of
$\Gamma$. The set $\bigtriangledown_Q(a)\rightleftharpoons
\bigcup\limits_{n\in\omega}Q^n(a,\Gamma)$  (respectively
$\bigtriangleup_Q(a)\rightleftharpoons
\bigcup\limits_{n\in\omega}Q^n(\Gamma,a)$) is called an {\em
upper} ({\em lower}) {\em
$Q$-cone}\index{$Q$-cone!upper}\index{$Q$-cone!lower}\index{$\bigtriangledown_Q(a)$}\index{$\bigtriangleup_Q(a)$}
of $a$. We call the $Q$-cones $\bigtriangledown_Q(a)$ and
$\bigtriangleup_Q(a)$ by {\em cones}\index{Cone} and denote by
$\bigtriangledown(a)$ and
$\bigtriangleup(a)$\index{$\bigtriangledown(a)$}\index{$\bigtriangleup(a)$}
respectively if $Q$ is fixed.

Recall \cite{SuLP, Su071, SuGP} that a countable acyclic directed
graph $\Gamma=\langle X;Q\rangle$ is said to be
\emph{powerful}\index{Digraph!powerful} if the following
conditions hold:

(a) the automorphism group of $\Gamma$ is {\em
transitive}\index{Group of automorphisms!transitive}, that is any
two vertices are connected by an automorphism;

(b) the formula $Q(x,y)$ is equivalent in the theory ${\rm
Th}(\Gamma)$ to a~disjunction of principal formulas;

(c) ${\rm acl}(\{a\})\cap\bigtriangleup_Q(a)=\{a\}$ for each
vertex $a\in X$;

(d) $\Gamma\models\forall x,y\:\exists z\:(Q(z,x)\wedge Q(z,y))$
(the \emph{pairwise intersection property}\index{Property!pairwise
intersection}).

\medskip
Below we define the property of powerfulness for the directed
graph $\Gamma$ in terms of the groupoid $\mathfrak{P}_{\nu(p)}$
for the unique $1$-type $p$ of the theory $T={\rm Th}(\Gamma)$
assuming that the theory is small.

At first we note that $U^-=\varnothing$ in view of Corollary 1.3
and so $\mathfrak{P}_{\nu(p)}$ is a monoid.

Since the formula $Q(x,y)$ is equivalent to some disjunction
$\bigvee\limits_{i=1}^n\theta_{u_i}(x,y)$, the acyclicity of
$\Gamma$ means that $0\notin u_{i_1}u_{i_2}\ldots u_{i_k}$ for any
$u_{i_1},\ldots,u_{i_k}\in\{u_1,\ldots,u_n\}$. The condition ${\rm
acl}(\{a\})\cap\bigtriangleup(a)=\{a\}$ is equivalent to that no
set $u^{-1}_{i_1}u^{-1}_{i_2}\ldots u^{-1}_{i_k}$ does not contain
almost deterministic elements. The pairwise intersection property
means that for any $u_i$, $i=1,\ldots,n$, and any $v\in U$ the set
$u_iv$ contains an element $u_j$. In particular, if $n=1$ then
$u_1\in u_1v$ for any $v\in U$. In this case we say that the
element $u_1$ {\em induces the pairwise intersection
property}\index{Element!inducing the pairwise intersection
property} or is a {\em ${\rm PIP}$-element}\index{${\rm
PIP}$-element}.

The characterizations above imply the following

\medskip
{\bf Proposition 7.1.} {\em A small theory $T$ of language
$\{Q^{(2)}\}$ is a theory of a powerful graph $\Gamma=\langle
X;Q\rangle$ if and only if $T$ has the unique $1$-type $p$ with a
regular labelling function $\nu(p)$ such that for some elements
$u_1,\ldots,u_n\in\rho_{\nu(p)}$ the following conditions are
satisfied:

$(1)$ $\vdash
Q(x,y)\leftrightarrow\bigvee\limits_{i=1}^n\theta_{u_i}(x,y)$;

$(2)$ $0\notin u_{i_1}u_{i_2}\ldots u_{i_k}$ for any
$u_{i_1},\ldots,u_{i_k}\in\{u_1,\ldots,u_n\}$;

$(3)$ for any $u_i$, $i=1,\ldots,n$, and any $v\in U$ the set
$u_iv$ contains an element $u_j$.}

\medskip
{\bf Definition.} A monoid $\mathfrak{P}_{\nu(p)}$ is called {\em
special}\index{Monoid!special} if $\rho_{\nu(p)}\cap
U^-\ne\varnothing$ and for any elements
$u_1,u_2,\ldots,u_n,v\in\rho_{\nu(p)}$, where
$u_1<0,\ldots,u_n<0$, $v\geq 0$, and for any element $u'\in
u_1u_2\ldots u_nv$ there is an element $v'\geq 0$ such that $u'\in
v'u_1u_2\ldots u_n$.

A special monoid $\mathfrak{P}_{\nu(p)}$ is called {\em ${\rm
PIP}$-special}\index{Monoid!${\rm PIP}$-special} if each negative
element $u\in\rho_{\nu(p)}$ is a  ${\rm PIP}$-element, i.~e.,
$u\in uv$ for any $v\in\rho_{\nu(p)}$.

\medskip
Having a special monoid (for a special small theory $T$) the
process of construction of a limit model over a type $p$ is
reduced to a sequence of $\theta_{u_n}$-extensions, $u_n<0$,
$n\in\omega$, of prime models over realizations of $p$: for any
limit model $\mathcal{M}$ over $p$ there is an elementary chain
$(\mathcal{M}(a_n))_{n\in\omega}$, $\models p(\bar{a}_n)$, such
that its union forms $\mathcal{M}$ and
$\models\theta_{u_n}(a_{n+1},a_n)$ is satisfied, $n\in\omega$. In
this case the isomorphism type of $\mathcal{M}$ is defined by the
sequence $(u_n)_{n\in\omega}$.

As shown in \cite{SuLP}, if a ${\rm PIP}$-special monoid exists
then, by adding of multiplace predicates, each prime model over a
tuple of realizations of $p$ is transformed to a model isomorphic
to $\mathcal{M}_p$. Thus, the type $p$ is connected with the
unique, up to isomorphism, prime model over realizations of $p$
and with some (finite, countable, or continual) number of limit
models over $p$, which is defined by some quotient for the set of
sequences $(u_n)_{n\in\omega}$, $u_n\in U^-\cap\rho_{\nu(p)}$,
$n\in\omega$. The action of these quotients is defined by some
identifications $(w\approx w')$ of words in the alphabet
$U^-\cap\rho_{\nu(p)}$ such that if $w=u_1\ldots u_m$ and
$w'=u'_1\ldots u'_n$ then for any $v\in U^{\geq
0}\cap\rho_{\nu(p)}$ and $u_0\in u_1\ldots u_mv$ there exists
$v'\in U^{\geq 0}\cap\rho_{\nu(p)}$ with $u_0\in v'u'_1u'_2\ldots
u'_n$.

To conclude this section we describe some connections of
$I_{\nu(p)}$-monoids with the strict order property.

\medskip
{\bf Definition.} Let $T$ be a theory with a type $p$ having the
model $\mathcal{M}_p$, $\mathfrak{P}_{\nu(p)}$ be an
$I_{\nu(p)}$-groupoid, and $X$ be a subset of $\rho_{\nu(p)}$
having a cardinality $\lambda$. We say that $X$ is (formula) {\em
definable}\index{Set!formula definable} if for a realization $a$
of $p$ the set of solutions of $L_{\lambda^+,\omega}$-formula
$\varphi(a,y)\rightleftharpoons\bigvee\limits_{u\in
X}\theta_{u}(a,y)$ in $\mathcal{M}_p$ is
$L_{\omega,\omega}$-definable in $\mathcal{M}_p$ by a formula
$\psi(a,y)$. In this case we say that the formula $\psi(x,y)$ {\em
witnesses}\index{Formula!witnessing on definability of set} on the
definability of $X$.

We say that a groupoid $\mathfrak{P}_{\nu(p)}$ {\em generates the
strict order property}\index{Groupoid!generating the strict order
property} if, for some definable set $X\subseteq\rho_{\nu(p)}$,
for a witnessing formula $\varphi(x,y)$, and for some realizations
$a$ and $b$ of $p$ satisfying $\models\theta_v(b,a)$ with a label
$v\in\rho_{\nu(p)}$, the inclusion
$\varphi(a,\mathcal{M}_p)\subset\varphi(b,\mathcal{M}_p)$ holds.

\medskip
{\bf Proposition 7.2.} {\em If $T$ is a small theory with a type
$p$, and the groupoid $\mathfrak{P}_{\nu(p)}$ has a definable set
$X\subseteq\rho_{\nu(p)}$ containing an element $u<0$ with $u\cdot
X \subseteq X$, then $\mathfrak{P}_{\nu(p)}$ generates the strict
order property.}

\medskip
{\em Proof.} Take a definable set $Y=X\cup\{0\}$ and consider a
witnessing formula $\varphi(x,y)$. Since $u\cdot X \subseteq X$
then $u\cdot Y \subseteq Y$ and, for any realizations $a$ and $b$
of $p$ with $\mathcal{M}\models\theta_{u}(b,a)$, we have
$\varphi(a,\mathcal{M}_p)\subseteq\varphi(b,\mathcal{M}_p)$. At
the same time, $0\in Y$ implies $b\in\varphi(b,\mathcal{M}_p)$,
and if $b\in\varphi(a,\mathcal{M}_p)$ then $a$ isolates $b$ that
is impossible by $u<0$. Thus,
$\varphi(a,\mathcal{M}_p)\subset\varphi(b,\mathcal{M}_p)$ and
$\mathfrak{P}_{\nu(p)}$ generates the strict order
property.~$\Box$

\medskip
{\bf Corollary 7.3.} {\em Let $T$ be a small theory with a type
$p$, and for some nonempty finite set $X\subseteq
U^-\cap\rho_{\nu(p)}$ there be a natural number $n$ such that
$X^{n+1}\subseteq\bigcup\limits_{i=1}^n X^i$,  $X^1=X$,
$X^{i+1}=X^i\cdot X$. Then the groupoid $\mathfrak{P}_{\nu(p)}$
generates the strict order property.}

\medskip
{\em Proof.} Clearly, the finite set $X$ is definable and the sets
$X^i$ and $Y\rightleftharpoons\bigcup\limits_{i=1}^n X^i$ are also
definable. Since $X^{n+1}\subseteq Y$ then for any element $u\in
X$ we have $u\cdot Y\subseteq Y$. Since $u<0$ then, by Proposition
7.2, the groupoid $\mathfrak{P}_{\nu(p)}$ generates the strict
order property.~$\Box$

\medskip
{\bf Corollary 7.4.} {\em If $T$ is a small theory with a type $p$
and $U^-\cap\rho_{\nu(p)}$ is a nonempty finite set then the
groupoid $\mathfrak{P}_{\nu(p)}$ generates the strict order
property.}

\medskip
{\em Proof.} Consider the set $X=U^-\cap\rho_{\nu(p)}$. As $X$ is
finite it is definable. Since $X$ contains all negative labels in
$\rho_{\nu(p)}$, by Proposition 1.4, we have $u\cdot X\subseteq X$
for any $u<0$ in $\rho_{\nu(p)}$. Therefore, by Proposition 7.2,
the groupoid $\mathfrak{P}_{\nu(p)}$ generates the strict order
property.~$\Box$

\medskip
\section{Partial groupoid  of binary isolating \\ formulas  on
a set of realizations \\ of a family of $1$-types of a complete
theory}
\medskip

In this section, the results above for a structure of a type are
generalized for a structure on a set of realizations for a family
of types.

Let $R$ be a nonempty family of types in $S^1(T)$. We denote by
$\nu(R)$\index{$\nu(R)$} a regular family of labelling
functions\index{$\rho_{\nu(R)}$}
$$\nu(p,q)\mbox{\rm : }{\rm PF}(p,q)/{\rm PE}(p,q)\to U,\,\,\,p,q\in
R,$$ $$\rho_{\nu(R)}\rightleftharpoons\bigcup\limits_{p,q\in
R}\rho_{\nu(p,q)}.$$

Similarly Proposition 3.1, we obtain that, having atomic models
$\mathcal{M}_p$ for all types $p\in R$ (for instance, if $T$ is
small), the function $P$, being partial for $|R|>1$, on the set
$R\times(\mathcal{P}(U)\setminus\{\varnothing\})\times R$, which
maps each tuple of triples
$(p_1,u_1,p_2),\ldots,(p_k,u_k,p_{k+1})$, where
$u_1\in\rho_{\nu(p_1,p_2)},\ldots,u_k\in\rho_{\nu(p_k,p_{k+1})}$,
to the set of triples $(p_1,v,p_{k+1})$, where $v\in
P(p_1,u_1,p_2,u_2,\ldots,p_k,u_k,p_{k+1})$, is {\em left
semi-associative}:\index{Function!semi-associative}
\begin{equation}
\begin{array}{c}P(P(p_1,u_1,p_2,u_2,p_3),u_3,p_4)=P(p_1,u_1,p_2,u_2,p_3,u_3,p_4)\supseteq
\\
\supseteq P(p_1,u_1,P(p_2,u_2,p_3,u_3,p_4))
\end{array}
\end{equation}
for $u_1\in\rho_{\nu(p_1,p_2)}$, $u_2\in\rho_{\nu(p_2,p_3)}$,
$u_3\in\rho_{\nu(p_3,p_4)}$.

Having the models $\mathcal{M}_p$ we consider the semi-associative
structure $\mathfrak{P}_{\nu(R)}\rightleftharpoons\langle
R\times(\mathcal{P}(U)\setminus\{\varnothing\})\times
R;\,\cdot\rangle$\index{$\mathfrak{P}_{\nu(R)}$} with the partial
operation $\cdot$ such that
$$(p_1,X_1,p_2)\cdot(p_2,X_2,p_3)=\bigcup\{(p_1,u_1,p_2)\cdot(p_2,u_2,p_3)\mid u_1\in X_1,u_2\in X_2\},$$
$$
(p_1,u_1,p_2)\cdot(p_2,u_2,p_3)=\{(p_1,v,p_3)\mid v\in
P(p_1,u_1,p_2,u_2,p_3)\},$$
$$u_1\in\rho_{\nu(p_1,p_2)},
u_2\in\rho_{\nu(p_2,p_3)}.
$$
The groupoids $\mathfrak{P}_{\nu(p)}$, $p\in R$, are naturally
embeddable in this structure. The structure
$\mathfrak{P}_{\nu(R)}$ is called a {\em join of
groupoids}\index{Join of!groupoids} $\mathfrak{P}_{\nu(p)}$, $p\in
R$, relative to the family $\nu(R)$ of labelling functions and it
is denoted by
$\bigoplus\limits_{\nu(R)}\mathfrak{P}_{\nu(p)}$\index{$\bigoplus\limits_{\nu(R)}\mathfrak{P}_{\nu(p)}$}.
If $\rho_{\nu(p,q)}=\varnothing$ for all $p\ne q$ the join
$\bigoplus\limits_{\nu(R)}\mathfrak{P}_{\nu(p)}$ is {\em
free},\index{Join!of groupoids!free}, it is isomorphically
represented as the disjoint union of the groupoids
$\mathfrak{P}_{\nu(p)}$ and denoted by $\bigsqcup\limits_{p\in
R}\mathfrak{P}_{\nu(p)}$.\index{$\bigsqcup\limits_{p\in
R}\mathfrak{P}_{\nu(p)}$}

By $(8)$, we have

\medskip
{\bf Proposition 8.1.} {\em For any complete theory $T$, for any
nonempty family $R\subset S(T)$ of $1$-types having models
$\mathcal{M}_p$ for each $p\in P$, and for any regular family
$\nu(R)$ of labelling functions, each $n$-ary partial operation
$$P(p_1,\cdot,p_2,\cdot,p_3\ldots,p_{n},\cdot,p_{n+1})$$ on the set
$\mathcal{P}(U)\setminus\{\varnothing\}$ is interpretable by a
term of the structure
$\bigoplus\limits_{p\in\nu(R)}\mathfrak{P}_{\nu(p)}$ with fixed
types $p_1,\ldots,p_{n+1}\in R$.}

\medskip
By Proposition 1.4, we obtain the following analogue of
Proposition 3.3.

\medskip
{\bf Proposition 8.2.} {\em For any complete theory $T$, for any
nonempty family $R\subset S(T)$ of $1$-types, and for any regular
family $\nu(R)$ of labelling functions, the restriction of the
structure $\mathfrak P_{\nu(R)}$ to the set of negative {\rm
(}respectively non-positive, non-negative{\rm )} labels is closed
under the partial operation $\cdot$.}

\medskip
In view of Proposition 8.2, the structure $\mathfrak P_{\nu(R)}$
has substructures $\mathfrak P^{\leq 0}_{\nu(R)}$\index{$\mathfrak
P^{\leq 0}_{\nu(R)}$} and $\mathfrak P^{\geq
0}_{\nu(R)}$,\index{$\mathfrak P^{\geq 0}_{\nu(R)}$} generated by
triples $(p,u,q)$ with $u\leq 0$ and $u\geq 0$ respectively,
$p,q\in R$. Here, for any triple $(p,u,q)$ in $\mathfrak P^{\geq
0}_{\nu(R)}$ the triple $(q,u^{-1},p)$ is also attributed to
$\mathfrak P^{\geq 0}_{\nu(R)}$.

A structure $\mathfrak P_{\nu(R)}$ is called ({\em almost}) {\em
deterministic}\index{Structure!deterministic}\index{Structure!almost
deterministic} if the set $(p,u,q)\cdot(q,v,r)$ is a singleton
(finite) for any triples $(p,u,q)$ and $(q,v,r)$ in $\mathfrak
P_{\nu(R)}$ with $u\in\rho_{\nu(p,q)}$ and $v\in\rho_{\nu(q,r)}$.

The deterministic structure $\mathfrak P_{\nu(R)}$ is generated by
the structure $\mathfrak P'_{\nu(R)}=\langle R\times U\times R;\,
\odot\rangle$,\index{$\mathfrak P'_{\nu(R)}$} where
$(p,u,q)\cdot(q,v,r)=\{(p,u,q)\odot(q,v,r)\}$ for $p,q,r\in R$,
$u,v\in U$.

Adapting the proof of Proposition 4.1 to a family $R$ of $1$-types
we obtain

\medskip
{\bf Proposition 8.3.} {\em For any complete theory $T$, for any
nonempty family $R\subset S(T)$ of $1$-types having models
$\mathcal{M}_p$ for each $p\in P$, and for any regular family
$\nu(R)$ of labelling functions, the following conditions are
equivalent:

$(1)$ the relation $I_R$ is transitive for any model
$\mathcal{M}\models T$;

$(2)$ the structure $\mathfrak P_{\nu(R)}$ is almost
deterministic.}

\medskip
Note that the absence of principal edges linking distinct
realizations of types in $R$ is equivalent to the antisymmetry of
the relation $I_R$. Since $I_R$ reflexive (by the formula
$(x\approx y)$), the definition of the family $\nu(R)$ and
Propositions 1.4, 8.3 imply

\medskip
{\bf Corollary 8.4.} {\em For any complete theory $T$, for any
nonempty family $R\subset S(T)$ of $1$-types having models
$\mathcal{M}_p$ for each $p\in P$, and for any regular family
$\nu(R)$ of labelling functions, the following conditions are
equivalent:

$(1)$ the relation $I_R$ is a partial order on the set of
realizations of types of $R$ in any model $\mathcal{M}\models T$;

$(2)$ the structure $\mathfrak P_{\nu(R)}$ is almost deterministic
and $\rho_{\nu(R)}\subseteq U^{\leq 0}$.

The partial order $I_R$ is identical if and only if
$\rho_{\nu(R)}=\{0\}$. The non-identical partial order $I_R$ has
infinite chains if and only if $|\rho_{\nu(p)}|>1$ for some $p\in
R$ or there is a sequence  $p_n$, $n\in\omega$, of pairwise
distinct types in $R$ such that $|\rho_{\nu(p_n,p_{n+1})}|\geq 1$,
$n\in\omega$, or $|\rho_{\nu(p_{n+1},p_{n})}|\geq 1$,
$n\in\omega$.}

\medskip
Lemma 1.1 and Proposition 8.3 imply

\medskip
{\bf Corollary 8.5.} {\em For any complete theory $T$, for any
nonempty family $R\subset S(T)$ of $1$-types having models
$\mathcal{M}_p$ for each $p\in P$, and for any regular family
$\nu(R)$ of labelling functions, the following conditions are
equivalent:

$(1)$ $I_R$ is an equivalence relation on the set of realizations
of types of $R$ in any model $\mathcal{M}\models T$;

$(2)$ the structure $\mathfrak P_{\nu(R)}$ is almost deterministic
and $\rho_{\nu(p)}\subseteq U^{\geq 0}$.}

\medskip
An element $u\in U$ is called {\em {\rm (}almost{\em )}
deterministic}\index{Element!deterministic}\index{Element!almost
deterministic} with respect to the regular family $\nu(R)$ of
labelling functions if, for some realization $a$ of a type in $R$
and for some type $q\in R$, the formula $\theta_{({\rm
tp}(a),u,q)}(a,y)$ is consistent and has a unique solution (has
finitely many solutions).

Repeating the proof of Proposition 4.7 we have

\medskip
{\bf Proposition 8.6.} {\em For any structure $\mathfrak
P_{\nu(R)}$ its restriction $\mathfrak
P_{\nu(R),d}$\index{$\mathfrak P_{\nu(R),d}$} {\rm (}respectively
$\mathfrak P_{\nu(R),{\rm ad}}${\rm )}\index{$\mathfrak
P_{\nu(R),{\rm ad}}$} to the set of {\rm (}almost{\rm )}
deterministic elements is closed under the partial operation of
the structure $\mathfrak P_{\nu(R)}$.}

\medskip
Using the proof of Proposition 4.9 the following proposition
holds.

\medskip
{\bf Proposition 8.7.} {\em If for the types $p,q\in S^1(T)$ the
models $\mathcal{M}_p$ and $\mathcal{M}_q$ exist then an element
$u\geq 0$ in $\rho_{\nu(p,q)}$ is deterministic if and only if
$(q,u^{-1},p)\cdot (p,u,q)=\{(q,0,q)\}$.}

\medskip
{\bf Proposition 8.8.} {\em If the structure $\mathfrak
P_{\nu(R)}$ is deterministic then the structure $\mathfrak
P'_{\nu(R)}$ is a join of groups if and only if each set
$\rho_{\nu(p)}$, $p\in R$, consists of non-negative elements.}

\medskip
{\em Proof} repeats the proof of Proposition 4.10 for each set
$\rho_{\nu(p)}$.~$\Box$

\medskip
{\bf Corollary 8.9.} {\em If $R$ is a nonempty family of $1$-types
in $S^1(T)$, there are models $\mathcal{M}_p$ for $p\in R$,
$\mathfrak P_{\nu(R)}$ is a deterministic structure, and
$\mathfrak P'_{\nu(R)}$ is a join of groups, then all elements in
$\mathfrak P'_{\nu(p)}$, $p\in R$, are deterministic.}

\medskip
{\em Proof.} Since, by Proposition 8.8, the sets $\rho_{\nu(p)}$
consist of non-negative elements, the determinacy of the structure
$\mathfrak P_{\nu(R)}$ and Proposition 8.7 imply that each element
in $\mathfrak P'_{\nu(p)}$, $p\in R$, is deterministic.~$\Box$

\medskip
Repeating the proof of Proposition 4.12 we obtain

\medskip
{\bf Proposition 8.10.} {\em If $R$ is a nonempty family of
$1$-types in $S^1(T)$, there exists models $\mathcal{M}_p$ for
$p\in R$, and $\nu(R)$ is a regular family of labelling functions,
then for the structure $\mathfrak P_{\nu(R)}$ the set
$\rho_{\nu(R),d}^{\geq 0}$ of all non-negative deterministic
elements $u$ in $\rho_{\nu(R)}$, for which the elements $u^{-1}$
are also deterministic, forms the deterministic substructure
$\mathfrak G^{\geq 0}_{\nu(R),d}$\index{$\mathfrak G^{\geq
0}_{\nu(R),d}$} of $\mathfrak P_{\nu(R)}$ such that $(\mathfrak
G^{\geq 0}_{\nu(R),d})'$ is a join of groups.}

\medskip
The results above substantiate the transformation of the diagram
in Figure 1 replacing the type $p$ by a nonempty family
$R\subseteq S^1(\varnothing)$.

\medskip
\section{$I_\mathcal{R}$-structures}
\medskip

{\bf Definition.} Let $\mathcal{R}$ be a nonempty set,
$$U=U^-\,\dot{\cup}\,\{0\}\,\dot{\cup}\,U^+$$
be an alphabet consisting of a set $U^-$\index{$U^-$} of {\em
negative elements}\index{Element!negative}, of a set
$U^+$\index{$U^+$} of {\em positive
elements}\index{Element!positive} and a zero $0$. If $p$ and $q$
are elements in $\mathcal{R}$, we write $u<0$ and $(p,u,q)<0$ for
any $u\in U^-$, $u>0$ and $(p,u,q)>0$ for any $u\in U^+$. For the
set $\mathcal{R}^2$ of all pairs $(p,q)$, $p,q\in \mathcal{R}$, we
consider a {\em regular} family\index{Family!of sets!regular}
$\mu(\mathcal{R})$\index{$\mu(\mathcal{R})$} of sets
$\mu(p,q)\subseteq U$ such that

\medskip
${\small\bullet}$ $0\in\mu(p,q)$ if and only if $p=q$;

\medskip
${\small\bullet}$ $\mu(p,p)\cap\mu(q,q)=\{0\}$ for $p\ne q$;

\medskip
${\small\bullet}$ $\mu(p,q)\cap\mu(p',q')=\varnothing$ for $p\ne
q$ and $(p,q)\ne(p',q')$;

\medskip
${\small\bullet}$ $\bigcup\limits_{p,q\in\mathcal{R}}\mu(p,q)=U$.

\medskip
Further we write $\mu(p)$ instead of $\mu(p,p)$, and considering a
partial operation $\cdot$ on the set
$\mathcal{R}\times(\mathcal{P}(U)\setminus\{\varnothing\})\times\mathcal{R}$
we shall write, as above, $(p,u,q)\cdot(q,v,r)$ instead of
$(p,\{u\},q)\cdot(q,\{v\},r)$.

A left semi-associative structure $\mathfrak P=\langle
\mathcal{R}\times(\mathcal{P}(U)\setminus\{\varnothing\})\times
\mathcal{R};\,\cdot\rangle$ with a regular family
$\mu(\mathcal{R})$ of sets is called an {\em
$I_\mathcal{R}$-structure}\index{$I_\mathcal{R}$-structure} if the
partial operation $\cdot$ of $\mathfrak P$ has values
$(p,X,q)\cdot(p',Y,q')$ only for $p'=q$, $\varnothing\ne
X\subseteq\mu(p,q)$, $\varnothing\ne Y\subseteq\mu(p',q')$, and is
generated by the partial function $\cdot$ for elements in $U$
where $(p,x,q)\cdot(q,y,r)$ forms a nonempty set of triples
$(p,z,r)$, $z\in\mu(p,r)$, if $x\in\mu(p,q)$ and $y\in\mu(q,r)$:
for any sets $X,Y\in\mathcal{P}(U)\setminus\{\varnothing\}$,
$\varnothing\ne X\subseteq\mu(p,q)$, $\varnothing\ne
Y\subseteq\mu(q,r)$,
$$
(p,X,q)\cdot(q,Y,r)=\bigcup\{(p,x,q)\cdot(q,y,r)\mid x\in X,y\in
Y\},
$$
as well as the following conditions hold:

\medskip
${\small\bullet}$ each restriction $\mathfrak
P_{\mu(p)}$\index{$\mathfrak P_{\mu(p)}$} of $\mathfrak P$ to the
set
$\{p\}\times(\mathcal{P}(\mu(p))\setminus\{\varnothing\})\times\{p\}$
is isomorphic to an $I$-groupoid with the universe
$\mathcal{P}(\mu(p))\setminus\{\varnothing\}$, $p\in\mathcal{R}$;

\medskip
${\small\bullet}$ if $u\in\mu(p,q)$ and $u<0$ then the sets
$(p,u,q\cdot(q,v,r)$ and $(r,v',p)\cdot(p,u,q)$ consist of
negative elements for any $v\in \mu(q,r)$ and $v'\in(r,p)$;

\medskip
${\small\bullet}$ if $u\in\mu(p,q)$, $v\in\mu(q,r)$, $u>0$, and
$v>0$, then the set $(p,u,q)\cdot(q,v,r)$ consists of non-negative
elements;

\medskip
${\small\bullet}$ for any element $u\in\mu(p,q)$ with $u>0$ there
is the unique {\em inverse}\index{Element!inverse} element
$u^{-1}\in\mu(q,p)$, $u^{-1}>0$, such that $(p,0,p\in
(p,u,q)\cdot(q,u^{-1},p)$ and
$(q,0,q)\in(q,u^{-1},p)\cdot(p,u,q)$;

\medskip
${\small\bullet}$ if an element $(p,u,r)$ is positive and belongs
to the set $(p,v_1,q)\cdot(q,v_2,r)$ then the element
$(r,u^{-1},p)$ belongs to the set
$(r,v_2^{-1},q)\cdot(q,v_1^{-1},p)$;

\medskip
${\small\bullet}$ for any elements $(p,u_1,q),(q,u_2,r),(r,u_3,t)$
the following inclusion holds: $$((p,u_1,q)\cdot
(q,u_2,r))\cdot(r,u_3,t)\supseteq (p,u_1,q)\cdot((q,u_2,r)\cdot
(r,u_3,t)),$$ and the strict inclusion
$$((p,u_1,q)\cdot
(q,u_2,r))\cdot(r,u_3,t)\supset p,u_1,q)\cdot((q,u_2,r)\cdot
(r,u_3,t))$$ may be satisfied only for $u_1<0$ and
$|(q,u_2,r)\cdot(r,u_3,t)|\geq\omega$;

\medskip
${\small\bullet}$ the structure $\mathfrak P$ contains the {\em
deterministic}\index{Substructure!deterministic} substructure
$\mathfrak P^{\geq 0}_d$\index{$\mathfrak P^{\geq 0}_d$}, being
the restriction to the set
$$U^{\geq 0}_d=\{u\in
U^{\geq 0}\mid (q,u^{-1},p)\cdot(p,u,q)=\{(q,0,q)\}\mbox{ for some
}p,q\in\mathcal{R}\};$$ every set $(p,u,q)\cdot(q,v,r)$ is a
singleton for $u\in U^{\geq 0}_d\cap\mu(p,q)$ and $v\in U^{\geq
0}_d\cap\mu(q,r)$.

\medskip
By the definition, any $I_\mathcal{R}$-structure $\mathfrak P$
contains $I$-subgroupoids $\mathfrak P_{\mu(p)}$,
$p\in\mathcal{R}$, and $I_\mathcal{R}$-substructures $\mathfrak
P^{\leq 0}$\index{$\mathfrak P^{\leq 0}$} and $\mathfrak P^{\geq
0}$\index{$\mathfrak P^{\geq 0}$} being restrictions of $\mathfrak
P$ to the sets $U^{\leq 0}$ and $U^{\geq 0}$ respectively.

\medskip
{\bf Theorem 9.1.} {\em For any {\rm (}at most countable{\rm )}
$I_\mathcal{R}$-structure $\mathfrak P$ there exists a {\rm
(}small{\rm )} theory $T$ with a family $R\subset S(T)$ of
$1$-types and a regular family $\nu(R)$ of labelling functions
such that $\mathfrak P_{\nu(R)}=\mathfrak P$.}

\medskip
{\em Proof} follows the schema for the proof of Theorem 6.1
extended by the schema for the proof of Theorem 3.4.1 in
\cite{SuLP}. In view of bulkiness of this proof we only point out
the distinctive features leading to the proof of this theorem.

1. We introduce, for each symbol $p\in \mathcal{R}$, an unary
predicate $R_p$ intersecting with all predicates ${\rm Col}_n$,
$n\in\omega$, and forming, on the set of realizations of complete
$1$-type $p'(x)$, being isolated by the set
$\{R_p(x)\}\cup\{\neg{\rm Col}_n(x)\mid n\in\omega\}$, a structure
of isolating formulas correspondent to the $I$-groupoid $\mathfrak
P_{\mu(p)}$. Moreover, we suppose that predicates $R_p$ are
disjoint.

2. For the elements $u\in\mu(p,q)$ the predicates $Q_u$ link only
elements $a$ in $R_p$ with elements $b$ in $R_q$. Moreover, if
$u>0$ then ${\rm Col}(a)={\rm Col}(b)$, and if $u<0$ then ${\rm
Col}(a)\leq{\rm Col}(b)$ and the coloring ${\rm Col}$ is
$Q_u$-ordered.

3. The relation $Q^{\geq 0}=\bigcup\limits_{u\geq 0}Q_u$ is an
equivalence relation such that its classes are ordered by the
relation $Q^{<0}=\bigcup\limits_{u<0}Q_u$.~$\Box$

\medskip
In conclusion we note that, using the operation $\cdot^{\rm eq}$,
the constructions above can be transformed for an arbitrary family
of types in $S(T)$.

\end{document}